\documentclass[12pt,a4paper,twoside,final,notitlepage, leqno]{article}
\usepackage[english]{babel}
\usepackage[T1]{fontenc}  
\usepackage{graphicx}
\usepackage{amsmath,amsthm,epsfig,amsfonts,bbm}  

\newcommand{\pequationdeb}{$$ \left\{ \begin{minipage}[c]{130mm}}
\newcommand{\pequationfin}{\end{minipage}
                           \right. $$}

\def \smb {{\scriptstyle \bullet }}
\newcommand{\monitem}{ \smallskip \noindent $\bullet$ \quad  } 
\newcommand{\moneq}{\vspace*{-6pt} \begin{equation} \displaystyle } 
\newcommand{\moneqstar}{\vspace*{-6pt} \begin{equation*} \displaystyle } 
\newcommand{\monendstar}{\vspace*{-6pt} \end{equation*}   }
\newcommand{\monend}{\vspace*{-6pt} \end{equation}   }

\def\R{{\rm I}\! {\rm R}}

\def\section*#1{}

\def\resume{\if@twocolumn
\section*{R\'esum\'e}
\else \small
\quotation{\bf \it R\'esum\'e \rule[1mm]{1.5mm}{0.2mm}\vspace{0pt}}
\fi}
\def\endresume{\if@twocolumn\else\endquotation\fi}
\def\abstract{\if@twocolumn
\noindent\section*{{\bf Abstract}}
\else \small
\quotation{\noindent \bf {Abstract.} \rule[1mm]{1.5mm}{0.2mm}\vspace{0pt}}
\fi}
\def\endabstract{\if@twocolumn\else\endquotation\fi}

\usepackage{fancyhdr}
\renewcommand{\headrulewidth}{0pt}

\begin{document} 

\fancypagestyle{plain}{ \fancyfoot{} \renewcommand{\footrulewidth}{0pt}}
\fancypagestyle{plain}{ \fancyhead{} \renewcommand{\headrulewidth}{0pt}}

~

\bigskip \bigskip   \bigskip  
\bigskip \bigskip \bigskip  
 
\centerline {\bf \LARGE On hyperbolic systems with entropy velocity } 

 \bigskip  \centerline {\bf \LARGE  covariant under  the action of  a group }

\bigskip \bigskip \bigskip \bigskip  

\centerline { \large    Fran\c{c}ois Dubois $^{a,b}$ }   

\bigskip  

\centerline { \it  \small   
$^a$   Conservatoire National des Arts et M\'etiers,  }

\centerline { \it  \small   Laboratoire de M\'ecanique des Structures
 et des Syst\`emes Coupl\'es,  Paris, France. } 

\centerline { \it  \small   $^b$  Department of Mathematics, University  Paris-Sud,}

\centerline { \it  \small B\^at. 425, F-91405 Orsay Cedex, France } 

 \bigskip

\centerline { \it  \small  francois.dubois@math.u-psud.fr}

\smallskip 
\bigskip  \bigskip \bigskip

\centerline { {\rm  31 august 2014} 
\footnote {\rm  \small $\,\,$ 
This contribution  in honor of Jean-Marie Souriau  (1922 - 2012) is 
issued from two lectures  entitled ``Sur les syst\`emes hyperboliques
\`a vitesse entropique invariants sous l'action d'un groupe'' 
given on  tuesday 28 august  2012 
and ``On hyperbolic systems with entropy velocity  covariant
 under  the action of  a group'' 
given on  wednesday 21 august  2013  at the 56th and 57th
``Colloque International de Th\'eories Variationnelles'', 
La Baume, Aix en Provence, France. 
Published in the {\it Journal of Hyperbolic Differential Equations}, 
volume~12, number 4, pages~763-785, 2015,
doi: 10.1142/S0219891615500228. Edition 23 january 2016. }}

\bigskip  \bigskip \bigskip

\noindent  {\bf Abstract. } \qquad 
For hyperbolic systems of conservation laws in one space dimension with a mathematical
entropy, we define the notion of entropy velocity.
Then we give sufficient conditions for such a system to be covariant under the action of 
a group of space-time transformations. These conditions naturally introduce a
representation  of the group in the space of states.
We  construct such  hyperbolic system from the knowledge of data on
the manifold of null velocity. We apply these ideas for Galileo, Lorentz and circular
groups. We focus on particular non trivial examples for two
by  two systems of conservation laws. 

\bigskip \noindent 
   {\bf Mots-cl\'es} : Galileo group, Lorentz group, circular group, mathematical entropy.



 \newpage

\fancyhead[EC]{\sc{Fran\c{c}ois Dubois}} 
\fancyhead[OC]{\sc{On hyperbolic systems with entropy velocity... }} 
\fancyfoot[C]{\oldstylenums{\thepage}}

\bigskip \bigskip \noindent {\bf \large 1) \quad  Introduction}  

\fancyfoot[C]{\oldstylenums{\thepage}}
 

\monitem 
The first link between covariance under the action of a group
and  conservation laws is the classic result of  
Emmy Noether \cite{No18}: 
when the action is invariant under the action of a group, 
an associated moment is conserved along the time evolution. 
This result is well understood for systems that possess
some action. 
This result is also fundamental  for all the modern physics. 
It has been intensively  studied and generalized by Jean-Marie  Souriau 
in his  fundamental book \cite{So70}
(see also the work of Kostant \cite{Ko70}). 
The fact to consider Galilean relativity is 
 natural in this kind of work. 
In order to explicit the mass as an invariant, the 
Bargmann group \cite{Ba54} is necessary in this framework
and we refer to the    work of  de  Saxc\'e and Vall\'ee
 \cite{SV12}  for the introduction of  the
second principle of thermodynamics. 
Observe also 
that very early Souriau has considered 
not only Galilean relativity but also 
Einstein's relativity  \cite{So64, So77, So78}. 
%
An other contribution is due to Vall\'ee \cite {Va87} relative to continuum
media   in the framework of special relativity.

\monitem 
Our scientific background  is first concerned with 
 the approximation of gas dynamics problems 
with numerical methods    \cite{DD05}. The notion of first order hyperbolic 
system of  conservation laws is fundamental for this study and the
contribution of Peter lax (see {\it e.g.} \cite {La06})  is
essential for all that follows. 
The introduction of the second principle of thermodynamics in 
this framework is due to Godunov \cite{Go61} and to 
Friedrichs and Lax  \cite{FL71} with the notion of
``mathematical entropy''. A  fundamental result of existence
and uniqueness of an entropy
solution in the scalar case is due to Kru$\breve {\rm z} $kov \cite {Kr70}. 
Galilean conservation for  hyperbolic systems has been 
proposed by Ruggeri \cite{Ru89}. 
We have considered the same question with 
enforcing the structure of the mathematical entropy in 
\cite{Du01}.  We   consider an extension of this 
approach in this contribution. 

\monitem 
The summary of our study is the following. We 
 introduce the entropy velocity at Section~2
for hyperbolic systems with a mathematical entropy. 
Then we define in Section~3 
the notion of covariance for an hyperbolic system with entropy velocity 
under the action of a group 
and propose sufficient conditions that make in evidence 
 an  algebraic structure. 
The  manifold of null velocity  defined at 
Section~4 is a natural notion in ths context. 
Then precise constraints for the entropy variables
are derived at Section~5. In Section~6, we study the particular case
of   systems of order two. Some words of conclusion are proposed at Section~7.


\bigskip \bigskip   \noindent {\bf \large 2) \quad  Entropy velocity } 

\monitem
Let $N$ be a positive integer and $\, \Omega \,$ a  non void  convex open part 
of $ \, \R^N .$ 
The physical flux $ \, f(W) \in \R^N  \,$ of a state 
$\, W \in \Omega \,$ is a  regular function 
$ \, f : \, \Omega \longrightarrow \R^N . \,$ 
We study in the following the system of conservation laws in one space dimension 
with unknown 
$\, \R \times [0, \, +\infty[ \, \,  \ni (x, \, t)  \, \longmapsto \, W(x, \, t) \in \Omega \, $ 
associated with this nonlinear flux function: 
\moneq   \label{edp}
{{\partial W}\over{\partial t}} \,+\, {{\partial}\over{\partial x}}  f(W) \,=\, 0  \,. 
\monend  

\monitem
As proposed by  Godunov \cite{Go61}, Friedrichs and Lax \cite{FL71}
and Boillat \cite{Bo74}, 
a mathematical entropy for the system (\ref{edp}) is a regular strictly convex
function  $ \, \eta : \, \Omega \longrightarrow \R \,$ 
such that for any regular solution $\, W(x, \, t) \,$ of the conservation law
(\ref{edp}), there exists a complementary conservation equation 
\moneq   \label{cons-entropy}
{{\partial \eta}\over{\partial t}} \,+\, {{\partial \zeta}\over{\partial x}}  
 \,=\, 0  \,. 
\monend  
The  space derivative of the second term in  (\ref{cons-entropy}) is 
the entropy flux, a regular function $ \, \zeta :   \Omega \longrightarrow \R \,$ 
that satisfies  the identity 
\moneq   \label{flux-entropy}
{\rm d} \zeta \, \equiv \, {\rm d} \eta \, \smb \, {\rm d} f  \,. 
\monend  
It is well known (see {\it e.g.} Godlewski - Raviart \cite{GR96} or \cite{DD05})
that a system of conservation law which  admits a mathematical  entropy is hyperbolic: the
jacobian $ \,  {\rm d} f(W) \,$ is diagonalisable in $ \, \R^N .\,$ 
Observe also that the restriction to only one space dimension
is not restrictive in principle if we adopt a point of view founded
on mathematical entropy. The generalization to several space dimensions
does not add in priciple new fundamental difficulties. 
As well known, the equality  (\ref{cons-entropy}) has to be replaced by an 
inequality (taken in an appropriate weak sense) when discontinuous  solutions
of the conservation law  (\ref{edp}) are considered. 
The co-vector  $ \, \varphi \,$  of entropy variables is the 
jacobian of the entropy relatively to the variation of the state 
$ \, W $:
\moneq   \label{var-entropic}
{\rm d} \eta \, \equiv \, \varphi \,\smb \,  {\rm d} W \,,   \qquad \forall \, W \in \Omega  \,. 
\monend  
The dual entropy $ \, \eta^* \,$ is the 
Legendre \cite{Le87} - Fenchel \cite{Fe49} -  Moreau \cite{Mo62}  
transform of the convex function $ \, \eta $. 
It is a function of the entropy variables:
\moneq   \label{dual-entropy}
\eta^* (\varphi)  \,\,\, \,\,\, \equiv \,\,\,\,\,\, 
{\rm sup}_{_{_{_{_{    \!\!\!\!\!\!\!\!\!\!\!\!\!\!   \normalsize  W \in \Omega}}}}} \,  
\big( \varphi \, \smb\, W - \eta(W) \big) \, 
\monend  
and we have 
\moneq   \label{dual-entropy-grad}
{\rm d} \eta^* (\varphi)  \, = \, {\rm d} \varphi \, \smb\, W \, . 
\monend

\monitem
In the following,  
we define the ``entropy velocity'' $\, u \,$ as the quotient of the entropy flux
divided by the entropy. In other terms, we suppose that the system  
 (\ref{edp})  (\ref{cons-entropy}) admits an entropy velocity $ \, u \,$
if there exists a regular function 
 $ \, u : \, \Omega \longrightarrow \R \,$ such that 
the entropy flux $ \, \zeta \, $ can be written under the form 
 $ \, \eta \, u  $:
\moneq   \label{def-vitesse}
\zeta(W) \, \equiv \, \eta(W) \,\, u(W) \,, \qquad \forall \, W \in \Omega  \,. 
\monend  
This definition has been previously proposed in \cite{Du01}. For gas dynamics 
this entropy velocity is the usual velocity because the specific entropy
is advected by the flow (see  {\it e.g.}  \cite{DD05} or \cite{GR96}). 
Once the velocity is defined, it is natural, as suggested 
by  Ruggeri \cite{Ru89},  to decompose  the physical  flux $ \, f \, $ into a convective part 
$ \, u \, W \,$ and a complementary  contribution  $ \, g $. 
We define the  ``thermodynamic flux''  $ \, g(W) \,  $ by the relation 
\moneq   \label{def-flux-thermo}
f(W) \, \equiv \,  u(W) \, W \,+\, g(W) \,, \qquad \forall \, W \in \Omega  \,. 
\monend  

\newpage 
 
\smallskip \monitem   { \bf Proposition 1.  Compatibility relation  } 

\noindent 
With the above framework, the compatibility relation  (\ref{flux-entropy}) 
is satisfied  if and only if 
\moneq   \label{compati}
\varphi \, \smb \, {\rm d}g \,+ \, \eta^* \,  {\rm d}u \, \equiv \, 0 \, .   
\monend  
%

\smallskip \noindent      {\bf Proof of Proposition 1. } 

\noindent 
The proof has been explicited in \cite{Du01};  we detail it for completeness. 
We have, due to the relations (\ref{flux-entropy}) to  (\ref{def-flux-thermo}), 
the following calculus:

\noindent $ \displaystyle 
{\rm d} \zeta \, -\, {\rm d} \eta \, \smb \, {\rm d} f  \, = \, 
{\rm d} ( \eta \, u ) \, - \,  {\rm d} \eta \,\smb \,  {\rm d}  \big(  u \, W \,+\, g  \big)   \, = \, 
\eta \, {\rm d} u  \, + \, u \,  {\rm d} \eta \, - \, 
\varphi \, \smb \, \big[ 
( {\rm d} u )  \, W  \, + \, u \,  {\rm d}  W   \, + \,  {\rm d} g  \, \big] 
$ 

\noindent $ \displaystyle \qquad    = \,\, 
\big( \eta \, - \, \varphi \, \smb \, W \big) \,  {\rm d} u  
\, + \, u \, \big(  {\rm d} \eta  \,-\,  \varphi \, \smb \,  {\rm d} W \big)
 \, - \,  \varphi \, \smb \,  {\rm d} g   \,\,\, = \, \,
- \eta^* \,  {\rm d} u    \, - \,  \varphi \, \smb \,  {\rm d} g $ 

\smallskip \noindent 
and the conclusion is clear.  $ \hfill \square $  

\bigskip \bigskip   \noindent {\bf \large 3) \quad  Covariance under the action of a group } 

\monitem
We introduce a group  $ \, {\cal G} \, $ with one real parameter $ \, \theta \, $  
composed by transformations $ \,  G_\theta \, $ of space-time  
 defined by a two by two matrix:
\moneq   \label{groupe-general}
G_\theta \, = \,  \begin{pmatrix}  \alpha_\theta &     \delta_\theta \cr 
\beta_\theta &   \gamma_\theta \end{pmatrix}  , \qquad  G_\theta \in  {\cal G}
 \,, \quad  \theta \in \R \, . \monend   
At a fixed $ \, \theta $, the transformation 
$ \, (x, \, t) \longmapsto (x', \, t') \,$ is defined by the relations:
\moneq   \label{transfo-theta}  
 \begin{pmatrix}  x' \cr t'  \end{pmatrix} \, = \, 
G_\theta \,  \begin{pmatrix}  x \cr t  \end{pmatrix} \, ; \qquad  
 \left\{ \begin{array} [c]{l} \displaystyle 
x' \, = \,  \alpha_\theta \, x \, + \,   \delta_\theta \, t  
 \\ \displaystyle    
t' \,\, = \,  \beta_\theta \, x \, + \,   \gamma_\theta \, t  \, . 
\end{array} \right. \monend  
In the following, we   consider three particular groups, 
 defined by the conditions
\moneqstar   
 \alpha_\theta \,=\,  \gamma_\theta  \,=\, C_\theta \,, \quad 
 \beta_\theta \,=\, {{\varepsilon}\over{c}} \, S_\theta \,, \quad 
 \delta_\theta  \,=\, c \,  S_\theta \, 
\monendstar 
{\it id est} by the following matrices
\moneq   \label{groupe-particulier}
G_\theta \, = \,  \begin{pmatrix}  C_\theta &   c\, S_\theta  \cr 
\displaystyle  {{\varepsilon}\over{c}} \, S_\theta &    C_\theta \end{pmatrix}  \, . \monend   
The Galileo group corresponds to 
\moneq   \label{gpe-galilee}
\varepsilon \, = \, 0 \,, \qquad C_\theta \, \equiv \, 1 
\,, \qquad S_\theta \, \equiv \, \theta  \,, 
\monend   
the Lorentz group to 
\moneq   \label{gpe-lorentz}
\varepsilon \, = \, 1 \,, \qquad C_\theta \, \equiv \, {\rm cosh} \, \theta  
\,, \qquad S_\theta \, \equiv \,  {\rm sinh}  \, \theta  \,, 
\monend   
and a third one here called the ``circular group'' to 
\moneq   \label{gpe-circular}
\varepsilon \, = \, -1 \,, \qquad C_\theta \, \equiv \, {\rm cos}  \, \theta  
\,, \qquad S_\theta \, \equiv \,  {\rm sin}  \, \theta  \,.  
\monend   
We remark that we have the elementary relations, 
\moneq   \label{trigo}
C_\theta ^2 \, - \, \varepsilon  \, S_\theta ^2 \, \equiv \, 1 \,, \quad 
{{{\rm d} S_\theta}\over{{\rm d} \theta}} \, \equiv \,  C_\theta \,, \quad 
{{{\rm d} C_\theta}\over{{\rm d} \theta}} \, \equiv \, \varepsilon  \, S_\theta \,,  \quad 
\varepsilon \in \{ -1 ,\, 0 , \, 1 \} \, .   
\monend   

\monitem
We say here that the system (\ref{edp}) associated with the   conservation of entropy
\moneq   \label{cons-entropy-vitesse}
{{\partial \eta}\over{\partial t}} \,+\, {{\partial}\over{\partial x}}  
\big( \eta \, u \big)  \,=\, 0  \, 
\monend  
is covariant under the action of the group $ \, \cal G \, $ 
if for each $ \, \theta \in \R \, $ and for each state $ \, W \in \Omega ,$ 
 we can define a new state $ \, W' \,$ 
such that after the change of variables 
(\ref{transfo-theta}), the system of equations  (\ref{edp})  (\ref{cons-entropy-vitesse}) 
take the form 
\moneq   \label{edp-prime}
 \left\{ \begin{array} [c]{l} \displaystyle 
\, \, \, \,   {{\partial W'}\over{\partial t'}} \, \quad + \quad 
\, {{\partial}\over{\partial x'}}  \,   f(W') \quad \,=\, 0   
  \\ \displaystyle   \vspace{-.5cm}  ~  \\ \displaystyle
{{\partial}\over{\partial t'}} \,  \eta(W') \,+\, {{\partial}\over{\partial x'}}  
 \,  \big( \eta(W') \, u(W') \big)  \,=\, 0  \,  
\end{array} \right. \monend  
with the {\it same}  flux function $ \, f(\smb) , \,$  the   {\it same}  mathematical entropy
 $\, \eta(\smb) \, $ and the   {\it same} entropy velocity  $\, u (\smb) . \, $

\monitem
We present now an  elementary calculus. We have the chain rule 
\moneqstar 
 \left\{ \begin{array} [c]{l} \displaystyle 
{{\partial}\over{\partial x}} \,\,=\,\, 
{{\partial x'}\over{\partial x}} \, {{\partial}\over{\partial x'}} \,+ \,  
{{\partial t'}\over{\partial x}} \, {{\partial}\over{\partial t'}}  \,\,=\,\, 
\alpha_\theta \, {{\partial}\over{\partial x'}} \,+ \,  
\beta_\theta  \, {{\partial}\over{\partial t'}}
  \\ \displaystyle   \vspace{-.3cm}  ~   \\ \displaystyle 
{{\partial}\over{\partial t}} \, \,\,=\,\, 
{{\partial x'}\over{\partial t}} \, {{\partial}\over{\partial x'}} \,+ \,  
{{\partial t'}\over{\partial t}} \, {{\partial}\over{\partial t'}}  \,\,=\,\, 
 \delta_\theta  \, {{\partial}\over{\partial x'}} \,+ \,  
\gamma_\theta \, {{\partial}\over{\partial t'}} \, .
\end{array} \right. \monendstar   
We inject these operators inside the relation  (\ref{edp}):

\smallskip \noindent $ \displaystyle \qquad \qquad
\Big(  \delta_\theta  \, {{\partial}\over{\partial x'}} \,+ \,  
\gamma_\theta \, {{\partial}\over{\partial t'}} \Big) \, W 
 \,+ \, \Big( \alpha_\theta \, {{\partial}\over{\partial x'}} \,+ \,  
  \beta_\theta \, {{\partial}\over{\partial t'}}  \Big) \, f(W)  
\,=\, 0  $

\smallskip \noindent and we put in evidence the new partial derivatives: 

\smallskip \noindent $ \displaystyle \qquad \qquad
 {{\partial}\over{\partial t'}}  \Big( \gamma_\theta \, W \, + \, 
  \beta_\theta \, f(W) \Big) \, + \, 
 {{\partial}\over{\partial x'}}  \Big(  \delta_\theta  \,  W \, + \, 
 \alpha_\theta \, f(W) \Big)  \,=\, 0  \, . $

\smallskip \noindent
We introduce a bijective linear operator $ \, Y_\theta \,$   
depending on $ \, \theta \in \R , \,$
acting from $ \, \R^N \,$ into  $ \, \R^N   \,$  and independent on 
$\, x' \,$ and $ \, t' $:
\moneq   \label{Y-theta}
Y_\theta \in {\cal GL} (  \R^N )  \,, \qquad \theta \in \R \, . 
\monend  
After applying $ \, Y_\theta \,$ to  the previous equation, we obtain 
a new conservation law equivalent to  (\ref{edp}):
\moneq   \label{edp-prime-vec}
 {{\partial}\over{\partial t'}}  \Big( \gamma_\theta \,\,   Y_\theta \,\smb \,  W \, + \, 
 \beta_\theta  \,\,  Y_\theta \,\smb \, f(W) \Big) \, + \, 
 {{\partial}\over{\partial x'}}  \Big(  \delta_\theta  \,\,   Y_\theta \,\smb \, W \, + \, 
 \alpha_\theta \,\,  Y_\theta \,\smb \, f(W) \Big)  \,=\, 0  \, .  
\monend    
A {\it sufficient} condition to establish the first equation of 
(\ref{edp-prime}) is to identify the arguments of the two partial
derivatives $ \, \partial_{t'} \,$ and  $ \, \partial_{x'} \,  $ in the previous equation: 
\moneq   \label{transfo-vecteurs}
 \left\{ \begin{array} [c]{l} \displaystyle 
 W'  \quad \,\,\, =\,   \gamma_\theta \, \,   Y_\theta \,\smb \,  W \, + \, 
   \beta_\theta   \,\,   Y_\theta \,\smb \, f(W)  
  \\ \displaystyle
f(W') \,=\,     \delta_\theta \,\,    Y_\theta \,\smb \, W \, + \, 
 \alpha_\theta \,\,   Y_\theta \,\smb \, f(W)   \, .  
\end{array} \right. \monend  
The relation (\ref{transfo-vecteurs}) can be written in a more compact form: 
\moneq   \label{transfo-vecteurs-2} 
\begin{pmatrix}  f(W')  \cr W'  \end{pmatrix} 
\,=\, G_\theta \,\,   Y_\theta \,\smb \, 
\begin{pmatrix}  f(W)  \cr W  \end{pmatrix} 
\, \equiv \,  G_\theta \,\,  \begin{pmatrix}   
 Y_\theta \,\smb \,  f(W)  \cr  Y_\theta \,\smb \, W  \end{pmatrix}  \, . 
 \monend  
Under the ``linearity  hypothesis'' (\ref{Y-theta}), the relation 
 (\ref{transfo-vecteurs-2}) gives   geometrical constraints 
for the   transformation 
$ \, \Omega \ni W \longmapsto W' \in \Omega  \,$ 
associated with the action of the group $ \, \cal G $.

\monitem
We make an analogous calculus for the second equation of 
(\ref{edp-prime})  relative to entropy. After the transformation  of  partial derivatives, 
the relation  (\ref{cons-entropy-vitesse}) takes the form

\smallskip \noindent $ \displaystyle \qquad \qquad
\Big(  \delta_\theta  \, {{\partial}\over{\partial x'}} \,+ \,  
\gamma_\theta \, {{\partial}\over{\partial t'}} \Big) \, \eta(W) 
 \,+ \, \Big( \alpha_\theta \, {{\partial}\over{\partial x'}} \,+ \,  
  \beta_\theta \, {{\partial}\over{\partial t'}}  \Big) \, 
\big( \eta(W) \, u(W) \big) \,=\, 0  $

\smallskip \noindent
and  can be written as
\moneq   \label{edp-prime-entrop} 
 {{\partial}\over{\partial t'}}  \Big( \gamma_\theta \,\,   \eta(W)  \, + \, 
  \beta_\theta \, \,    \eta(W) \, u(W)   \Big) \, + \, 
 {{\partial}\over{\partial x'}}  \Big(  \delta_\theta  \, \,   \eta(W)  \, + \, 
 \alpha_\theta \, \,  \eta(W) \, u(W)   \Big)  \,=\, 0  \, . 
 \monend  
This is a scalar conservation law. This equation is identical to the second equation 
of (\ref{edp-prime}) 
if the following {\it sufficient} conditions are satisfied:  
\moneq   \label{transfo-entropie}
 \left\{ \begin{array} [c]{l} \displaystyle 
 \quad  \,\, \,   \eta(W')    \quad \,\,\, =\,   \gamma_\theta \, \,    \eta(W) \,   \, + \, 
   \beta_\theta   \,\,     \eta(W) \, u(W)
  \\ \displaystyle
     \eta(W') \, u(W')  \,=\,     \delta_\theta \,\,    \eta(W)  \, + \, 
 \alpha_\theta \,\,   \eta(W) \, u(W) \, .  
\end{array} \right. \monend  
With a notation that makes the action of the group $ \, \cal G \, $ explicit:
\moneq   \label{transfo-entropie-2} 
\begin{pmatrix}    \eta(W') \, u(W')  \cr   \eta(W')   \end{pmatrix} 
\,=\, G_\theta \,\,   
\begin{pmatrix}    \eta(W) \, u(W) \cr  \eta(W)   \end{pmatrix}    \, . 
 \monend  
%

\smallskip  \monitem   { \bf Proposition 2.  Linear representation  } 

\noindent 
If the relation (\ref{transfo-vecteurs-2})  is satisfied for every $ \, \theta \in \R, \,$ 
then the application $ \,  \, {\cal G }   \ni  G_\theta \longmapsto \, Y_\theta 
\in  {\cal GL } (\R^N) \,    $ is a linear representation of the group  $ \, {\cal G }$: 
\moneq   \label{rep-linear}
Y_{\theta +  \theta'} \,=\, Y_{\theta} \,\smb\, Y_{\theta'}     \, , \qquad 
\theta , \, \theta' \in \R \, .  
 \monend  
%

\smallskip \noindent      {\bf Proof of Proposition 2. } 

\noindent 
We first consider the relation  (\ref{transfo-vecteurs-2}) 
and  an analogous relation obtained by replacing the variable $ \, \theta \,$
by  the   variable $ \, \theta' $: 
\moneqstar 
\begin{pmatrix}  f(W'')  \cr W''  \end{pmatrix} 
\,=\, G_\theta'\,\,   Y_\theta' \,\smb \, 
\begin{pmatrix}  f(W')  \cr W'  \end{pmatrix} 
\, \equiv \,  G_\theta' \,\,  \begin{pmatrix}   
 Y_\theta' \,\smb \,  f(W')  \cr  Y_\theta' \,\smb \, W'  \end{pmatrix} \, .   
 \monendstar  
We compose this relation with   (\ref{transfo-vecteurs-2}) and we obtain
\moneqstar 
\begin{pmatrix}  f(W'')  \cr W''  \end{pmatrix} 
\,=\, G_{\theta'}\,\,\,  G_\theta\,\, 
\begin{pmatrix}  Y_{\theta'} \,\smb\, Y_{\theta} \, \smb \, f(W)  \cr 
 Y_{\theta'} \,\smb\, Y_{\theta}  \, \smb \, W  \end{pmatrix}  \, .   
 \monendstar  
When we apply the relation  (\ref{transfo-vecteurs-2}) with the argument 
$ \, \theta + \theta' $, we have directly
\moneqstar 
\begin{pmatrix}  f(W'')  \cr W''  \end{pmatrix} 
\,=\, G_{\theta + \theta'}\,\,\,  
\begin{pmatrix}  Y_{\theta + \theta'} \,\smb\,  f(W)  \cr 
   Y_{\theta + \theta'} \, \smb \, W  \end{pmatrix}  \, .   
 \monendstar  
We look precisely to the two previous relations. The scalars
$ \, \alpha_\theta $, $ \, \beta_\theta $,  $ \, \gamma_\theta \, $ and
 $ \, \delta_\theta \, $ commute with the linear operator 
$ \, Y_\theta \,$ and we obtain the relation (\ref{rep-linear}) 
because the two previous equalities are  satisfied for an arbitrary state $W$. 
  $ \hfill \square $  

\bigskip \bigskip   \noindent {\bf \large 4) \quad  Manifold of null velocity  } 

\monitem
In the following we denominate  by    $ \, \Omega_0 \,$
the set  of states with a  velocity equal to zero.
Because the mapping $ \, \Omega \ni W \longmapsto u(W) \in \R \,$ is scalar, it is natural 
to suppose that  $ \, \Omega_0 \,$ is a manifold in $\, \R^N \,$ of codimension one.  
We observe also that the flux function $ \, f(\smb) \,$ is reduced to its
thermodynamic contribution $ \, g(\smb) \,$ on the manifold of null velocity.
Our program is to construct the states $W$,  the flux function $ f(\smb) $
and the mathematical entropy  $ \eta(\smb) $
when  we suppose that the flux function 
$ \,  \Omega_0 \ni W_0 \longmapsto g_0(W_0) \in \R^N \, $ 
and the mathematical entropy 
$ \,  \Omega_0 \ni W_0 \longmapsto \sigma(W_0) \in \R \, $ 
are  given on the manifold of null
velocity.  
The idea is to use the group invariance to link a given state $ \, W \, $
with an appropriate state  $ \, W_0 \, $ on the manifold $ \, \Omega_0 \,$
as presented in Figure~1. The conditions (\ref{transfo-vecteurs-2}) and 
(\ref{transfo-entropie-2}) are geometric  constraints that will allow 
essentially to solve the problem.

 \bigskip    
\smallskip   \smallskip                   
\centerline { \includegraphics[width=.45 \textwidth] {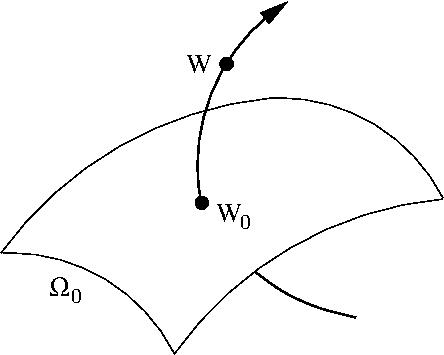} }   

\noindent  {\bf Figure 1}. \quad 
Fiber above the manifold of null velocity $\, \Omega_0 .$  Thanks to the action of 
the group, a state  $ \, W_0 \, $  of null velocity $\, \Omega_0 \,$ 
generates a current state $W$. 
\smallskip \smallskip 

\smallskip  \monitem   { \bf Proposition 3.  Velocity field } 

\noindent 
With the   hypothesis  (\ref{transfo-vecteurs}) and if the mathematical entropy
$ \, \sigma \,$ on the  manifold of null velocity is not the null function, 
the velocity field $ \, u(W) \,$ 
is necessarily given by the relation 
\moneq   \label{groupe-vitesse} 
   u(W) \,= \, {{\delta_\theta}\over{\gamma_\theta}} \, . 
 \monend  
%

\smallskip  \noindent      {\bf Proof of Proposition 3. } 

\noindent 
We consider a given state $ \, W_0 \,$ on the   manifold
$ \,  \Omega_0 \, $ and a running state $ \, W \in \Omega $ obtained by the relations 
  (\ref{transfo-vecteurs}). With the actual notations, we have 
\moneq  \label{transfo-vecteurs-3}
 W \, =\,   \gamma_\theta \, \,   Y_\theta \,\smb \,  W_0 \, + \, 
   \beta_\theta   \,\,   Y_\theta \,\smb \, g_0(W_0)  \,, \qquad  
f(W) \,=\,     \delta_\theta \,\,    Y_\theta \,\smb \, W_0 \, + \, 
 \alpha_\theta \,\,   Y_\theta \,\smb \, g_0(W_0)   \, .  
 \monend 
The analogous relations  (\ref{transfo-entropie}) relative to the entropy 
take a simple form  
\moneq  \label{transfo-entropie-3} 
 \eta(W)  \, =\,   \gamma_\theta \,     \sigma(W_0) \,\,, \qquad 
     \eta(W) \, u(W)  \,=\,     \delta_\theta \,\,    \sigma(W_0)   \, .   
\monend 
because the velocity is also null on $ \,  \Omega_0 $. 
We consider a state $W_0$ such that $ \, \sigma(W_0) \not=  0 .$ 
If $ \, \gamma_\theta \,$ is null, then $\, \eta(W) = 0 \, $ is null and 
$\, \delta_\theta = 0 \,$ because  $ \, \sigma(W_0) \not=  0 .$ 
This contradicts the fact that the matrix  
defining  $\, G_\theta \,$ by the relation  (\ref{groupe-general})
is invertible. Then $ \,  \gamma_\theta \not = 0 \,$ 
 and the relation 
(\ref{groupe-vitesse}) is   obtained by taking the ratio of the two  
relations of   (\ref{transfo-entropie-3}). 
  $ \hfill \square $  

\smallskip  \monitem 
With the example proposed in (\ref{groupe-particulier}), we have 
\moneq  \label{vitesse-groupe}
u \,=\, c \, {{S_\theta}\over{C_\theta}} \, .   
\monend 
We observe that we have also established the following property: 
 
\smallskip  \newpage \monitem   { \bf Proposition 4.  Entropy  field } 

\noindent With the above hypotheses, if the state $ \, W \,$ is given 
from  $ \, W_0 \in \Omega_0 \,$  
with the help  of the relations (\ref{transfo-vecteurs-3}), 
we have necessarily 
\moneq  \label{entropie} 
 \eta(W)  \, =\,   \gamma_\theta \,     \sigma(W_0)    \, .   
\monend 

\monitem { \bf Hypothesis 1.  
The   null-velocity manifold   $ \Omega_0 \, $ is  included in an hyperplane} 

\noindent 
We remark that if a convex restriction of the mathematical entropy  
$ \,  \Omega_0 \ni W_0 \longmapsto \sigma(W_0) \in \R \, $ 
is given on the manifold $ \,  \Omega_0 $, 
a natural hypothesis is to suppose  $ \,  \Omega_0 \, $
convex. We will do this hypothesis in the following. Moreover, 
we suppose  that $ \,  \Omega_0 \, $ is flat; this is expressed 
by our  hypothesis~1. 
 
%

\smallskip  \monitem   { \bf Proposition 5.  Thermodynamic flux } 
    
\noindent With the above framework (\ref{transfo-vecteurs-3}), 
the thermodynamic flux $ \, g(W) \,$ is given from the 
datum $ \, g_0(W_0) \,$ according to 
\moneq  \label{thermo-flux} 
 g (W)  \, = \, \big(  \alpha_\theta \, - \, u(W) \, \beta_\theta \,  \big)  \, 
Y_\theta \,\smb\,   g_0 (W_0)    \, .   
\monend 
%
 
\smallskip \noindent      {\bf Proof of Proposition 5. } 

\noindent 
The second relation of (\ref{transfo-vecteurs-3}) and the definition 
of the thermodynamic flux  (\ref{def-flux-thermo}) allows us to evaluate 
the function $\, g(\smb) .$ We have 

\smallskip  \noindent $ \displaystyle 
g(W)    \, = \, f(W) - u(W) \, W  $

\smallskip  \noindent $ \displaystyle \qquad \,\,\,\,  
 \, = \,   \delta_\theta \,\,    Y_\theta \,\smb \, W_0 \, + \, 
 \alpha_\theta \,\,   Y_\theta \,\smb \, g_0(W_0)   \, - \, 
u(W) \, \big(  \gamma_\theta \, \,   Y_\theta \,\smb \,  W_0 \, + \, 
   \beta_\theta   \,\,   Y_\theta \,\smb \, g_0(W_0) \big) \,  $ 

\smallskip  \noindent $ \displaystyle \qquad \,\,\,\,  
 \, = \,   \big(  \alpha_\theta \, - \, u(W) \, \beta_\theta \,  \big)  \, 
 Y_\theta \,\smb \, g_0(W_0) \,  $  \qquad due to (\ref{groupe-vitesse}) 

\smallskip  \noindent 
and the relation (\ref{thermo-flux}) is established. 
  $ \hfill \square $  

\monitem
For the groups proposed in (\ref{groupe-particulier}), we have 
$ \,\, \, \alpha_\theta \, - \, u(W) \, \beta_\theta \,=\, 
C_\theta - c \, {{S_\theta}\over{C_\theta}} \,  {{\varepsilon}\over{c}} \, S_\theta 
  \,=\,  {{1}\over{C_\theta}} $ \, \quad  due to (\ref{trigo}).  
Then the relation  (\ref{thermo-flux}) takes the form
\moneq  \label{thermo-flux-2} 
 g (W)  \, =\,   {{1}\over{C_\theta}}   \, Y_\theta \,\smb\,   g_0 (W_0)    \, .   
\monend

\bigskip \bigskip   \noindent {\bf \large 5) \quad  Entropy variables  } 

\monitem
Even if  the state and the entropy  
are only  partially known with the help of relations  (\ref{transfo-vecteurs-3}) 
and  (\ref{entropie}), its jacobian $\, \varphi \,$ introduced in  (\ref{var-entropic})
relative to the conserved variables $ \,  W \, $ can be essentially determined. 
We first consider the dynamical equations present in this study. 
First, the relation 
\moneq  \label{evol-group} 
 G_{\theta + \theta'}   \, =\,   G_{\theta} \, \smb \,   G_{\theta'}  
\monend 
and the relation (\ref{rep-linear}) for the matrices $ \, Y_\theta \,$ 
can be differentiated relatively to  $ \, \theta $. We have 
\moneq  \label{derive-groupe}  
{{{\rm d}G_\theta}\over{{\rm d}\theta}} \,=\, G'_0 \, G_\theta 
 \,=\, G_\theta \, G'_0   \,, \qquad  
{{{\rm d}Y_\theta}\over{{\rm d}\theta}} \,=\, Y'_0 \, Y_\theta 
 \,=\, Y_\theta \, Y'_0    \, .   
\monend  
With the choice  (\ref{groupe-particulier}) 
of one classical group, we have 
\moneq   \label{derivee-zero-groupe}
G'_0 \, = \,  \begin{pmatrix}  0  &   c  \cr 
\displaystyle  {{\varepsilon}\over{c}}   &  0  \end{pmatrix}  \,  \monend   
and we can precise how to extend the differential relatations 
(\ref{evol-group}) and (\ref{derive-groupe}).

\smallskip  \monitem   { \bf Proposition 6.  Differential relations for the entropy variables } 

\noindent With the above context, we  have 
\moneq   \label{derive-theta-W0}
G'_0 \, \begin{pmatrix}  u \, \eta^* + \varphi \, g \cr  \eta^*   \end{pmatrix}  
\, {\rm d}\theta  \,+\,  \begin{pmatrix} \varphi \, Y'_0 \, f(W) \cr 
 \varphi \, Y'_0 \,\,  W  \end{pmatrix}  \, {\rm d}\theta  \,+\, G_\theta \, 
 \begin{pmatrix}  \varphi \, Y_\theta \, {\rm d}g_0 \cr
 \varphi \, Y_\theta \, {\rm d}W_0 - {\rm d}\sigma  \end{pmatrix}  \,=\, 0 
\,  \monend   
%
 
\smallskip  \noindent      {\bf Proof of Proposition 6. } 

\noindent 
We write the relation (\ref{transfo-vecteurs-3}) under a matricial form 
\moneq   \label{transfo-vecteurs-4} 
\begin{pmatrix}  f(W) \cr W  \end{pmatrix}  \,=\, \, G_\theta \, 
\begin{pmatrix}  Y_\theta \, g_0(W_0)  \cr  Y_\theta \, W_0  \end{pmatrix}  \,.
 \monend   
We differentiate the relation  (\ref{transfo-vecteurs-4}), taking into account the 
two contributions  relative to  
  the group $\, {\cal G} \,$ on one hand and along the manifold 
$ \, \Omega_0 \,$ on the other hand,  as illustrated on the figure~1. Then 

\smallskip  \noindent $ \displaystyle 
\begin{pmatrix}  {\rm d}f  \cr  {\rm d} W  \end{pmatrix}   \,=\, 
 G_0'  \,\,  G_\theta \, 
\begin{pmatrix}   Y_\theta \, \, g_0  \cr  Y_\theta \, \, W_0   \end{pmatrix} \,  {\rm d}\theta 
\,+\,   G_\theta \, \begin{pmatrix}  Y'_0 \,\, Y_\theta \,\, g_0 \cr 
 Y'_0 \,\, Y_\theta \,\, W_0  \end{pmatrix} \,  {\rm d}\theta    \,+\,  G_\theta \, 
\begin{pmatrix}   Y_\theta \,\, {\rm d} g_0  \cr   Y_\theta \,\, {\rm d} W_0  \end{pmatrix}
 \,$ 

\smallskip  \noindent 
and because 

\smallskip  \noindent $ \displaystyle 
  G_\theta \, \begin{pmatrix}  Y'_0 \, \, Y_\theta  \,\, g_0 \cr 
 Y'_0 \,  \,Y_\theta  \,\, W_0  \end{pmatrix} \,=\, 
 \begin{pmatrix} \alpha_\theta \,  Y'_0  \,\, Y_\theta  \,\, g_0 
\,+\, u\, \,  \gamma_\theta  \,\,  Y'_0  \,\, Y_\theta  \,\, W_0 \cr
 \beta_\theta \ \,,  Y'_0  \,\, Y_\theta  \,\, g_0 \,+\, 
 \gamma_\theta  \,\,  Y'_0 \,\, Y_\theta  \,\, W_0  \end{pmatrix} \,=\, 
 \begin{pmatrix} Y'_0 \, \big(  \alpha_\theta \,\, Y_\theta \,\, g_0 
\,+\, u\, \, \gamma_\theta \,\,  Y_\theta \,\, W_0 \big) \cr 
 Y'_0 \, \big(  \beta_\theta \, \, Y_\theta \,\, g_0 \,+\, 
 \gamma_\theta \,\, Y_\theta \,\, W_0 \big)  \end{pmatrix} 
$

\smallskip  \noindent $ \displaystyle \qquad \qquad \qquad  \quad  \,=\, 
 \begin{pmatrix}  Y'_0 \,\, f(W) \cr  Y'_0 \,\, W    \end{pmatrix} \, $ ,  

\smallskip  \noindent 
we have  
\moneq   \label{derive-transfo-vecteurs-4} 
\begin{pmatrix}  {\rm d}f  \cr  {\rm d} W  \end{pmatrix}   \,=\,  
 G'_0 \, \begin{pmatrix} f  \cr W  \end{pmatrix} \, {\rm d}\theta   \,+\, 
 \begin{pmatrix} Y'_0 \,\, f  \cr  Y'_0 \,\, W  \end{pmatrix}  \, {\rm d}\theta   \,+\,  
 G_\theta \,  \begin{pmatrix}   Y_\theta \,\, {\rm d} g_0  \cr  
 Y_\theta \,\, {\rm d} W_0  \end{pmatrix} \, . 
\monend  
We multiply each line of the relation (\ref{derive-transfo-vecteurs-4})
by the line vector $ \, \varphi .$ We remark that 

\smallskip  \noindent $ \displaystyle 
\varphi \, G'_0 \, \begin{pmatrix}  {\rm d}f  \cr  {\rm d} W  \end{pmatrix}   \,=\,  
\varphi \,  \begin{pmatrix} \alpha'_0 & (u \, \gamma)'_0 \cr
 \beta'_0 &  \gamma'_0   \end{pmatrix} \, 
\begin{pmatrix}  {\rm d}f  \cr  {\rm d} W  \end{pmatrix}   \,=\, 
\varphi \,   \begin{pmatrix} \alpha'_0 \,\, f \,+\,  (u \, \gamma)'_0 \,\, W \cr 
 \beta'_0 \, \, f +  \gamma'_0  \,\, W   \end{pmatrix} \, $

\smallskip  \noindent $ \displaystyle \qquad \qquad \qquad  \,=\, 
   \begin{pmatrix} \alpha'_0 \,\, \varphi \,\, f \,+\,  (u \, \gamma)'_0 \,\, \varphi \,\,  W \cr 
 \beta'_0   \,\, \varphi \, \,  f +  \gamma'_0   \,\, \varphi  \,\,  W   \end{pmatrix}  \,=\, 
G'_0 \,   \begin{pmatrix}  \varphi  \,\,  f \cr  \varphi  \,\,  W  \end{pmatrix} \, $  
%
%

\smallskip  \noindent 
and in a similar way, 

\smallskip  \noindent $ \displaystyle 
\varphi \,\, G_\theta \,\, Y_\theta \,\,  \begin{pmatrix}  {\rm d}g_0  \cr 
 {\rm d} W_0 \end{pmatrix}   \,=\,   \varphi \,\,  \begin{pmatrix}  
\alpha_\theta \,\,  Y_\theta \,\,  {\rm d}g_0 \,+\,  (u \, \gamma)_\theta \, \,
 Y_\theta \,\,  {\rm d}W_0  \cr \beta_\theta \,\,  Y_\theta \,\,  {\rm d}g_0 \,+\, 
 \gamma_\theta \,\,  Y_\theta \,\,  {\rm d}W_0 \end{pmatrix}      $

\smallskip  \noindent $ \displaystyle \qquad \qquad \qquad \quad \  \,=\, 
   \begin{pmatrix} \alpha_\theta \,\, \varphi \,\,  Y_\theta \,\,  {\rm d}g_0 \,+\, 
  (u \, \gamma)_\theta \,\,  \varphi \,\,   Y_\theta \,\,  {\rm d}W_0  \cr
 \beta_\theta \,\, \varphi \,\,  Y_\theta \, \,  {\rm d}g_0 \,+\, 
 \gamma_\theta \,\,  \varphi \,\,   Y_\theta \,\,  {\rm d}W_0  \end{pmatrix}  \,=\, 
G_\theta \,     \begin{pmatrix}  \varphi\, \,  Y_\theta \,\,  {\rm d}g_0 \cr 
 \varphi \,  \, Y_\theta \,\,  {\rm d}W_0   \end{pmatrix}  \, .    $ 
%
%

\smallskip  \noindent 
Then 
\moneq   \label{entropic-derive-transfo-vecteurs-4} 
 \begin{pmatrix} \varphi \,  {\rm d}f \cr  \varphi \,  {\rm d}W  \end{pmatrix}   \,=\, 
 G_0' \,  \begin{pmatrix}  \varphi \, \, f \cr  \varphi  \,\, W  \end{pmatrix} \, {\rm d}\theta   
\,+\,   \begin{pmatrix}  \varphi  \,\, Y'_0  \,\, f \cr  \varphi \, Y'_0  \,\, W
 \end{pmatrix}   \, {\rm d}\theta \,+\, G_\theta \,   \begin{pmatrix}  \varphi  \,\, Y_\theta \, 
 \,  {\rm d}g_0 \cr   \varphi  \,\, Y_\theta  \, \,  {\rm d}W_0  \end{pmatrix} \, . 
\monend  
\monitem
We develop the same calculus for the entropy. We start with a matricial form 
of the relation (\ref{transfo-entropie-3}):
\moneq   \label{transfo-entropie-4} 
 \begin{pmatrix} \eta \, u \cr \eta  \end{pmatrix}   \,=\, 
G_\theta \,  \begin{pmatrix} 0  \cr \sigma  \end{pmatrix}  
\monend  
Then 
\moneq   \label{derive-transfo-entropie-4} 
 \begin{pmatrix} {\rm d} ( \eta \, u ) \cr   {\rm d}  \eta  \end{pmatrix}   \,=\, 
  G'_0 \,  \begin{pmatrix} \eta \, u \cr \eta  \end{pmatrix} \,  {\rm d} \theta  
\,+\, G_\theta \,  \begin{pmatrix} 0  \cr  {\rm d}  \sigma  \end{pmatrix}   \,. 
\monend  
We make the difference between the relations
(\ref{derive-transfo-entropie-4}) and (\ref{entropic-derive-transfo-vecteurs-4}). 
Due to   (\ref{flux-entropy}) and (\ref{var-entropic}) 
the left hand side is equal to zero! 
Then the right hand side can be written as 

\smallskip  \noindent \qquad  \qquad $ \displaystyle 
 G_0' \,  \begin{pmatrix}  \varphi \, f -  \eta \, u 
\cr  \varphi \, W  - \eta   \end{pmatrix} \, {\rm d}\theta  
\,+ \,   \begin{pmatrix}  \varphi \,\, Y'_0 \,\, f \cr  \varphi \,\, Y'_0 \,\, W
 \end{pmatrix}   \, {\rm d}\theta \,+\, G_\theta \,   \begin{pmatrix}  \varphi \,\, Y_\theta  \,
\, {\rm d}g_0 \cr   \varphi \,\, Y_\theta  \, \, {\rm d}W_0  \, - \,  {\rm d}  \sigma 
 \end{pmatrix} \,.  $

\smallskip  \noindent
The relation (\ref{derive-theta-W0}) is now easily obtained by the explicitation of the dual 
entropy introduced in (\ref{dual-entropy})  and the thermodynamic flux defined in 
(\ref{def-flux-thermo}).  
  $ \hfill \square $  

\smallskip  \monitem   { \bf Proposition 7. Constraints for the  entropy variables } 

\noindent 
With the above framework, if we suppose moreover that the derivative 
of the group $\, \cal G \,$ for $ \, \theta = 0 \,$ is given by the relation
(\ref{derivee-zero-groupe}), 
 we have the following relations for the entropy variables 
\moneq   \label{variables-entropiques-cn1}
c \, \eta^* \,+\, \varphi \,\, Y_0' \,\, f(W) \,=\, 0   
\monend  {\vspace*{-.7 cm}}
\moneq \label{variables-entropiques-cn2} 
{{\varepsilon}\over{c}} \, \big( u \,  \eta^*  \,+\, \varphi \, g(W) \big) 
 \,+\, \varphi \,\, Y_0' \,\, W  \,=\, 0 
\monend  {\vspace*{-.7 cm}}
\moneq   \label{variables-entropiques-cn3}
\varphi \, Y_\theta \,\, {\rm d}g_0  \,=\, 0 
\monend  {\vspace*{-.7 cm}}
\moneq    \label{variables-entropiques-cn4}
\varphi \,\, Y_\theta \,\, {\rm d}W_0  \, - \, {\rm d}\sigma  \,=\, 0 \, . 
\monend   
%

\smallskip  \noindent      {\bf Proof of Proposition 7. } 

\noindent 
With the matrix $ \, G_0' \,$ chosen  according to (\ref{derivee-zero-groupe}), 
the relation (\ref{derive-theta-W0}) becomes

\smallskip  \noindent \qquad  \qquad $ \displaystyle 
\begin{pmatrix}  c \, \eta^*  \cr 
{{\varepsilon}\over{c}} \, \big( u \,  \eta^*  \,+\, \varphi \, g(W) \big) 
 \end{pmatrix}  \, {\rm d}\theta  \,+\,  \begin{pmatrix} \varphi \, Y'_0 \, f(W) \cr 
 \varphi \, Y'_0 \,\,  W  \end{pmatrix}  \, {\rm d}\theta  \,+\, G_\theta \, 
 \begin{pmatrix}  \varphi \, Y_\theta \, {\rm d}g_0 \cr
 \varphi \, Y_\theta \, {\rm d}W_0 - {\rm d}\sigma  \end{pmatrix}  \,=\, 0 $  

\smallskip  \noindent
and the relations  (\ref{variables-entropiques-cn1}) to  (\ref{variables-entropiques-cn4}) 
are clear.  
  $ \hfill \square $  

\monitem
The relation  (\ref{variables-entropiques-cn4}) can be expressed with the help of the
gradient $ \, \varphi_0 \,$ of the entropy restricted on the manifold $ \, \Omega_0 $: 
\moneq    \label{variables-entropiques-vit-nulles}
 {\rm d}\sigma \,\equiv \, \varphi_0 \, \smb \,  {\rm d}W_0 \, .  
\monend   
Then the relation  (\ref{variables-entropiques-cn4}) express that the 
difference $ \, \varphi \, Y_\theta - \varphi_0  \,$ belongs to the 
polar set $ \, \Omega_0^0 \,$ of  $ \, \Omega_0 $: we can write 
$ \,  \varphi \, Y_\theta - \varphi_0 \,=\, \mu \, r_0 \,$
for some scalar $\, \mu \,$ and some non-null linear form $ \, r_0 $
identically equal to zero on $ \, \Omega_0$. 
We make also a new hypothesis.

\monitem { \bf Hypothesis 2. The entropy $ \, \sigma \,$ and the thermodynamic flux 
$ \, g_0 \, $ are weakly  decoupled on the hyperplane $\Omega_0$}  

\noindent 
In other terms, 
\moneq    \label{bon-decouplage}
( \varphi_0 \,,\, {\rm d}g_0 \,\smb\, \rho_0 ) \,= \, 0 \,, \qquad \forall \rho_0 \in
\Omega_0 \, .  
\monend   
%

\smallskip   \monitem   { \bf Proposition 8. Explicitation of the  entropy variables } 

\noindent 
We make the  hypotheses 1 and 2 and suppose moreover  
\moneq    \label{pas-trivial}
\exists \, \rho_0 \in \Omega_0 \,, \quad  {\rm d}g_0 \,\smb\, \rho_0 \not= 0  \, .   
\monend   
Then the relations  (\ref{variables-entropiques-cn1}) to  (\ref{variables-entropiques-cn4}) 
imply a simple form for the entropy variables: 
\moneq    \label{var-entrop-simple} 
\varphi \, Y_\theta \,=\, \varphi_0   \, .   
\monend   
%

\smallskip  \noindent      {\bf Proof of Proposition 8. } 

\noindent 
Consider an arbitrary vector $ \, \rho_0 \in \Omega_0 $. 
When we apply the relation 
$ \,  \varphi \, Y_\theta - \varphi_0 \,=\, \mu \, r_0 \,$
established previously to the vector $ \,  {\rm d}g_0 \,\smb\, \rho_0 $,
we find zero due to the constraint  (\ref{variables-entropiques-cn3}) and 
the hypothesis~2 with (\ref{bon-decouplage}). 
Then $ \, \mu \, < r_0 \,,\,    {\rm d}g_0 \,\smb\, \rho_0 >  \, = \, 0 .\,$ 
 With the choice of $ \, \rho_0 \,$ suggested in (\ref{pas-trivial}), 
the duality product $ \,  < r_0 \,,\,    {\rm d}g_0 \,\smb\, \rho_0 >   \,$
is not null. So $ \, \mu = 0 \,$ and the relation  (\ref{var-entrop-simple})
is established. 
  $ \hfill \square $  

\monitem { \bf Hypothesis 3. The entropy $ \, \sigma \,$ and the thermodynamic flux 
$ \, g_0 \, $ are strongly  decoupled on the hyperplane $\Omega_0$}  

\noindent 
In other terms, 
\moneq    \label{fort-decouplage}
 \varphi_0 (W_0) \, \smb \, g_0 (\widetilde{W_0}) \,  \,= \,\, 0 \,, 
\qquad \forall \,   W_0, \, \widetilde{W_0}  \in \Omega_0 \, .  
\monend   
We observe that this strong hypothesis (\ref{fort-decouplage}) 
implies (\ref{bon-decouplage}) by differentiation relatively to 
the variable $ \, \widetilde{W_0} \,$ in $ \, \Omega_0 $. 
In consequence, we have a new property for the entropy variables. 

\smallskip   \monitem   { \bf Proposition 9. 
The entropy variables are in the polar set of the thermo flux  } 

\noindent 
We make the  hypotheses 1 and 3 and suppose also the technical 
condition (\ref {pas-trivial}) proposed  previously. 
Then we have a complete decoupling: 
\moneq    \label{varentrop-ortho-thermo-flux}
\varphi (W) \, \smb \, g({\widetilde{W}}) \, \equiv \, 0 \, 
\qquad \forall \,   W , \, \widetilde{W}  \in \Omega \, .   
\monend   
%
%

\smallskip  \noindent      {\bf Proof of Proposition 9. } 

\noindent 
We just have to precise the expression (\ref{thermo-flux-2}) of the thermodynamic flux. 
Then we have 

\smallskip  \noindent  \qquad $ \displaystyle 
\varphi (W) \, \smb \, g({\widetilde{W}}) \,= \, 
\varphi (W) \, \,Y_\theta \,\, Y_{-\theta} \,\, g({\widetilde{W}}) \,= \, 
\varphi_0 (W_0) \,\, {{1}\over{C_\theta}} \,  g_0({\widetilde{W_0}}) \,= \, 0
 $

\smallskip  \noindent
due to (\ref{varentrop-ortho-thermo-flux}).
  $ \hfill \square $  
 
\smallskip   \monitem   { \bf Proposition 10. Constraints for the representation $Y$ } 

\noindent 
Under the same  hypotheses  that the ones done for Proposition~9, we have 
\moneq   \label{variables-entropiques-cn1bis}
c \, \eta^* \,+\, C_\theta \,\, \varphi_0 \,\, Y_0' \,\, g_0  \,=\, 0   
\monend  {\vspace*{-.7 cm}}
\moneq \label{variables-entropiques-cn2bis} 
\varphi_0 \, \, Y_0' \,\, W_0  \,=\, 0 \,   \quad {\rm on}  \,  \, \Omega_0 \, . 
\monend 
%
%

\smallskip  \noindent      {\bf Proof of Proposition 10. } 

\noindent 
We first have the following calculus:

\smallskip  \noindent   \qquad $ \displaystyle 
\varphi \, \, Y_0' \, \, g 
\,=\, \varphi \, \, Y_\theta  \, \, Y_{-\theta} \,\, Y_0' \, \, g  
\,=\, \varphi \, \, Y_\theta  \, \, Y_0' \, \, Y_{-\theta} \,\,   g  
\,=\, \varphi_0  \, \, Y_0' \, \, {{1}\over{C_\theta}} \,  g_0   
\,=\,  {{1}\over{C_\theta}} \,  \varphi_0  \, \, Y_0' \, \,  g_0 \, . $

\smallskip  \noindent 
Then we multiply the equation (\ref{variables-entropiques-cn2}) by 
the opposite    of the velocity  $\, u \, $ and we add it to 
(\ref{variables-entropiques-cn1}). With the definition (\ref{def-flux-thermo})
of the thermodynamic flux, we have, taking into account the 
``trigonometrical relations''  (\ref{trigo}) and the expression  (\ref{groupe-vitesse})  of 
the velocity: 

\smallskip  \noindent  $ \displaystyle 
0 \,=\, 
\eta^* \, \Big( c - {{\varepsilon \, u^2}\over{c}} \Big) \,+ \, \varphi \, \, Y_0' \, \, g 
\,=\, \eta^* \,  \Big( c  - {{\varepsilon \, c \, S_\theta^2}\over{C_\theta^2}} \Big) 
 \,+ \,  {{1}\over{C_\theta}} \,  \varphi_0  \, \, Y_0' \, \,  g_0 
\,=\, 
  {{c\, \eta^*}\over{C_\theta^2}} \,+\, C_\theta \,  \varphi_0  \, \, Y_0' \, \,  g_0  $ 

\smallskip  \noindent 
and the relation  (\ref{variables-entropiques-cn1bis}) is established. 
We report now the expression  (\ref{transfo-vecteurs-3}) 
 of a current state $ \, W\, $ 
inside  the equation  (\ref{variables-entropiques-cn2}):

\smallskip  \noindent  $ \displaystyle 
0 \,=\,  {{\varepsilon}\over{c}} \, u \, \eta^* \,+\, 
\varphi \, \,   Y_0' \, \, Y_\theta   \,\smb \,\Big(  
 C_\theta \,  W_0 \, + \,   {{\varepsilon}\over{c}} \,  S_\theta \,  g_0 \Big) 
\,=\,  C_\theta \,\varphi \, \, Y_\theta  \, \, Y_0' \,\, W_0 
\, + \,  {{\varepsilon}\over{c}} \,\Big(  
u \, \eta^* \, +\, S_\theta \,
\varphi \, \, Y_\theta  \, \, Y_0' \,\, g_0 \Big)    $ 

\smallskip  \noindent  $ \displaystyle \quad 
= \,  C_\theta \,\varphi \, \, Y_\theta  \, \, Y_0' \,\, W_0 
\, + \,  {{\varepsilon}\over{c}} \,\Big(  
u \, \eta^* \, - \, S_\theta \,  {{c}\over{C_\theta}} \, \eta^*  \Big) 
\, = \,    C_\theta \,\,\varphi \, \, Y_\theta  \, \, Y_0' \,\, W_0  
\,+\, \varepsilon \, \eta^* \, \Big(   {{u}\over{c}} \,-\, 
 {{S_\theta}\over{C_\theta}}  \Big) $

\smallskip  \noindent  $ \displaystyle \quad 
= \,    C_\theta \,\,\varphi \, \, Y_\theta  \, \, Y_0' \,\, W_0    $ 

\smallskip  \noindent 
and the relation (\ref{variables-entropiques-cn2bis}) follows 
from (\ref{var-entrop-simple}). 
  $ \hfill \square $  

\smallskip   \monitem   { \bf Proposition 11. Dual entropy } 

\noindent 
Under the    hypotheses   proposed at  Proposition~9, we denote by 
$ \, \sigma^* \,$ the dual of the entropy $\, \sigma \,$ 
restricted on the manifold of null velocity:
\moneq   \label{dual-entropie-omega-zero}
\sigma^* \,\equiv \, \varphi_0 \,\smb \, W_0 \,-\, \sigma (W_0) \,  
\quad {\rm on}  \,  \, \Omega_0 \, , \qquad {\rm with } \quad \varphi_0 = \sigma'(W_0) \,. 
\monend  
Then we have 
\moneq   \label{dual-entropie}
\eta^* \,=\, C_\theta \, \, \sigma^*  \,   
\monend 
and the relation (\ref{variables-entropiques-cn1bis})
takes the form 
\moneq   \label{variables-entropiques-cn1ter}
c \, \,  \sigma^*   \,+\,  \, \varphi_0 \,\, Y_0' \,\, g_0  \,=\, 0    
  \qquad {\rm on}  \,  \, \Omega_0 \, .   
\monend  
%
%

\smallskip  \noindent      {\bf Proof of Proposition 11. } 

\noindent 
We have, thanks to (\ref{dual-entropy}),  (\ref{entropie}), (\ref{transfo-vecteurs-4}) 
and  (\ref{var-entrop-simple})~:

\smallskip  \noindent  \qquad $ \displaystyle 
\eta^* \, = \, \varphi \, \smb \, W \, - \, \eta(W) 
\,=\,  \varphi_0    \, \smb \, Y_{- \theta}  \, \, Y_{\theta} \, \big( C_\theta \, W_0 
\,+\, {{\varepsilon}\over{c}} \, S_\theta \, g_0  \big) 
 \, -  C_\theta \, \sigma  
\, \, =\,    C_\theta \,  \big(  \varphi_0   \, \smb \,  W_0 \,-\,  \sigma  \big)
$  
 
\smallskip  \noindent   due to  Hypothesis 3 (relation (\ref{fort-decouplage})). 
  The end of the proof is  clear. 
$ \hfill \square $  

\bigskip \bigskip   \noindent {\bf \large 6) \quad  Hyperbolic  systems of order two } 

\monitem
If we suppose $ \, N = 2 , \,$ the matrix $ \, Y_\theta \,$ is a two by two matrix. We
suppose here that this matrix has the same algebaic form that the form
(\ref{groupe-particulier}) suggested for $ \, G_\theta .\, $ We set for  this 
contribution 
\moneq   \label{matrice-Y-theta} 
Y_\theta \, = \,  \begin{pmatrix}  {\widetilde C}_\theta &    \displaystyle
{ {\widetilde \varepsilon}\over{a}} \,   {\widetilde S}_\theta   \cr 
\displaystyle  a \,  {\widetilde S}_\theta   &  $\,\,$  {\widetilde C}_\theta
\end{pmatrix} \,, \qquad \theta \in \R \,, \quad a > 0 \, .  \monend    
In the relation (\ref{matrice-Y-theta}), the functions $ \,  {\widetilde S}_\theta \,$
and  $ \, {\widetilde C}_\theta \,$ are analogous to the ones proposed 
in (\ref{gpe-galilee}),  (\ref{gpe-lorentz}),   (\ref{gpe-circular}) and 
 (\ref{trigo}) 
for  $ \, S_\theta \,$ and  $ \, C_\theta .\,$  We have in particular 
\moneq   \label{trigo-tilde} 
 {\widetilde C}_\theta^2 \, - \, {\widetilde \varepsilon} \, 
 {\widetilde S}_\theta^2 \,\, \equiv \,\, 1 \, . 
\monend    
The nilpotent representation of the group $\, {\cal G} \,$ corresponds to
$ \,  {\widetilde \varepsilon} = 0 ,\,$ the hyperbolic representation 
to $ \,  {\widetilde \varepsilon} = 1 \,$ and   the elliptic one 
to $ \,  {\widetilde \varepsilon} = -1 . \,$ The derivative $ \, Y_\theta ' \,$ of 
the matrix $ \, Y_\theta \,$  proposed in  (\ref{matrice-Y-theta}) at $ \,\theta
= 0 \,$ is given by 
\moneq   \label{matrice-Y-prime-zero} 
Y_0'  \, = \,  \begin{pmatrix}  0 &    \displaystyle { {\widetilde \varepsilon}\over{a}} \,    \cr 
\displaystyle  a \,     &  0  \end{pmatrix} \,.
 \monend    

\monitem
We particularize the system by a simple choice for the null-velocity manifold.
As in our previous work \cite{Du01}, we suppose 
\moneq   \label{variete-zero-Neq2} 
W_0 \, = \, \begin{pmatrix}  \rho_0  \cr 0   \end{pmatrix}  \, \, \in \Omega_0 \,   
 \monend    
and 
\moneq   \label{g-zero-Neq2} 
g_0 (W_0)  \, = \, \begin{pmatrix}  0  \cr p_0   \end{pmatrix}  \, 
\quad {\rm for} \,\, W_0  \in \Omega_0 \,  . 
\monend    
With this choice, the hypotheses 2 and 3 are satisfied. 
The variable $ \, p_0 \,$ in the right hand side of 
the relation (\ref{g-zero-Neq2}) is the thermodynamic pressure.

\smallskip   \monitem   { \bf Proposition 12. Pressure as a thermodynamic variable } 

\noindent With the above hypotheses, we have the relation
\moneq   \label{pression}
\sigma^* \,+\, {{\widetilde \varepsilon}\over{a \, c}} \,\,\sigma' \,\, p_0 \,=\, 0    \,  . 
 \monend    
%

\smallskip  \noindent      {\bf Proof of Proposition 12. } 

\noindent Just write the relation (\ref{variables-entropiques-cn1ter})
with the help of  (\ref{matrice-Y-prime-zero}) and   (\ref{variete-zero-Neq2})~:

\smallskip  \noindent  $ \displaystyle 
c \, \,  \sigma^*   \,+\,  \, \varphi_0 \,\, Y_0' \,\, g_0  \,=\, 
c \, \,  \sigma^*   \,+\,   \begin{pmatrix}  \sigma' & 0  \end{pmatrix} 
\,\, \begin{pmatrix}  0 &    \displaystyle { {\widetilde \varepsilon}\over{a}} \,    \cr 
\displaystyle  a \,     &  0  \end{pmatrix} \,\, 
  \begin{pmatrix}  0 \cr p_0  \end{pmatrix}  \,=\, 
c \, \,  \sigma^*   \,+\,   {{\widetilde \varepsilon}\over{a}} \, \sigma' \,\, p_0 $ 

\smallskip  \noindent  and the relation (\ref{pression}) is established.  
$ \hfill \square $  

\monitem  
 The nilpotent case  ($ \, \widetilde \varepsilon = 0 \,$) 
has no interest when $ \, N=2 \,$   because  the dual of the entropy is
necessary null, as we have 
remarked in  \cite{Du01} in the particular case of the Galileo group ($\varepsilon = 0$). 
We can explicit the first relations of (\ref{transfo-vecteurs-3}). With the notation
\moneq   \label{variables-conservatives}
W \, \equiv \,  \begin{pmatrix} \rho \cr J  \end{pmatrix}   \,  
\monend  
we have 
\moneq   \label{variables-conservatives-2}
\rho \, = \,  C_\theta \,   {\widetilde C}_\theta  \, \rho_0 
\,+\, {{\varepsilon \,  \widetilde \varepsilon}\over{a \, c}} \, 
 S_\theta \,   {\widetilde S}_\theta \, p_0    \,, \quad 
{{J}\over{a}}  \, = \,  C_\theta \,   {\widetilde S}_\theta  \, \rho_0
\,+\, {{\varepsilon }\over{a \, c}} \,  S_\theta \,   {\widetilde C}_\theta \, p_0    \,. 
\monend    

\monitem  
In the case of Galileo group, we have 
$\varepsilon = 0$, $  C_\theta = 1$,  $S_\theta = \theta$ and  $ u = c \, \theta$.
The explicitation of 
the parameter $ \, \theta \,$ and the state of null velocity $\, \rho_0 \,$
is   simple by resolution of the system (\ref{variables-conservatives-2}).
We have 
\moneqstar 
\rho \, = \,    {\widetilde C}_\theta  \, \rho_0    \,, \quad 
{{J}\over{a}}  \, = \,     {\widetilde S}_\theta  \, \rho_0  \,, \quad
\varepsilon = 0 \,, \quad {\rm Galileo} \,\, {\rm group}.   
\monendstar    
In the hyperbolic case, we have 
\moneqstar 
\tanh \theta \, = \, {{J}\over{a \, \rho}} \,, \quad
u \,=\, c \,  \theta \,, \quad \varepsilon = 0 \,, \,\, \widetilde \varepsilon = 1 \, 
\monendstar   
and in the elliptic one, 
\moneqstar 
{\rm tg }  \, \theta \, = \, {{J}\over{a \, \rho}} \,, \quad
u \,=\, c \,  \theta \,, \quad \varepsilon = 0 \,, \,\, \widetilde \varepsilon = -1 \, 
\monendstar    
as we have explained in   \cite{Du01}.

\monitem 
In the following, we focus on hyperbolic ($ \widetilde \varepsilon = 1 $)
and  elliptic  ($ \widetilde \varepsilon = -1 $) representations of the Lorentz 
 ($\varepsilon = 1$) and  circular ($\varepsilon = -1$) groups. With  other words, 
\moneq   \label{epsilon-tilde-carre}
 \varepsilon^2 \,\,= \,\, 
 \widetilde \varepsilon^2 \,\,= \,\, 1 \, .  
\monend  
We can replace the pressure 
$ \, p_0 \, $ with an appropriate thermodynamic quantity because 

\smallskip  \noindent \qquad \qquad  \qquad \qquad  \qquad \qquad  $ \displaystyle
{{\sigma^*}\over{\sigma'}} \,=\, {{\rho_0 \, \sigma' - \sigma}\over{\sigma'}}  \,=\,
\rho_0  \, - \,  {{\sigma}\over{\sigma'}}  \, .    $ 

\smallskip  \noindent 
Then the relations (\ref{variables-conservatives-2}) take the form 
\moneq   \label{variables-conservatives-3}
\rho \, = \,  \big( C_\theta \,   {\widetilde C}_\theta 
\,-\, \varepsilon \,  S_\theta \,   {\widetilde S}_\theta  \big) \, \rho_0 
\,+\, \varepsilon \,   S_\theta \,   {\widetilde S}_\theta \, {{\sigma}\over{\sigma'}} \, 
\,, \quad 
{{J}\over{a}}  \, = \,  \big(  C_\theta \,   {\widetilde S}_\theta  
\,-\, \varepsilon \,  \widetilde  \varepsilon\,  S_\theta \,   {\widetilde C}_\theta 
 \big) \,  \rho_0
\,+\,  \varepsilon  \, \widetilde  \varepsilon \, \,   S_\theta \,   {\widetilde C}_\theta \,
 {{\sigma}\over{\sigma'}} \,   \,. 
\monend    
For the Lorentz  ($\varepsilon = 1$)
or the circular  ($\varepsilon = -1$) group, the resolution of the system 
 (\ref{variables-conservatives-3}) with unknowns $\, \theta \,$ and 
$\, \rho_0 \, $ is absolutly non trivial. 
In order to be able to differentiate the solution  $\, (\theta ,\, 
 \rho_0 )\, $ relatively to $ \, (\rho \,,   J) ,\,$ we have the following result.

\smallskip    \monitem   { \bf Proposition 13. Jacobian } 

\noindent 
We denote by $ \, \zeta_0 \, $ the partial derivative of the ratio $ \, \sigma / \sigma' \,$
relatively to the  dentity $ \, \rho_0 \,$ at rest~:
\moneq   \label{zeta_zero}
\zeta_0  \,  \equiv \,  {{\partial}\over{\partial \rho_0}} 
\, \Big( {{\sigma}\over{\sigma'}} \Big) \, . 
 \monend    
If the hypothesis (\ref{epsilon-tilde-carre}) is satisfied and 
if the following jacobian determinant 
\moneq   \label{det-jaco}
\Delta  \, \equiv \, {{\partial \, ( \rho , \, J)}
\over {\partial  \,  (\theta ,\, \rho_0 )}} 
\,=\, a \,  {{\sigma}\over{\sigma'}} \, \Big[  (\varepsilon \,  \widetilde  \varepsilon 
-1) \,\, C_\theta^2 \,\,  {{\sigma^*}\over{\sigma}} \, \,
+\,\,  \widetilde  \varepsilon \, S_\theta^2 \, \zeta_0 
\,\, - \, \,(  C_\theta^2 \,+ \, \widetilde  \varepsilon \,  S_\theta^2 ) \Big] 
 \monend    
is not null, the solution  $\, (\theta ,\,  \rho_0 )\, $ of the system 
 (\ref{variables-conservatives-3}) is locally unique and we can differentiate
the parameters $ \, \theta \,$ and $ \, \rho_0 \,$ relatively
to the state $ \, W \, $ introduced in (\ref{variables-conservatives}).

\smallskip  \noindent      {\bf Proof of Proposition 13. } 

\noindent 
We have

\moneqstar    \left\{ \begin{array} [c]{rcl} \displaystyle 
  {{\partial \rho}\over{\partial \theta}} \,&=&\, \displaystyle 
(   \widetilde  \varepsilon \, - \,   \varepsilon ) \, 
  C_\theta  \,  \widetilde S_\theta  \, \,  \rho_0 
\,+\,  \varepsilon \big(   C_\theta  \,  \widetilde S_\theta  \, 
+ \,  S_\theta  \,  \widetilde C_\theta  \big) \, \,   {{\sigma}\over{\sigma'}} 
  \\ \displaystyle   \vspace{-.3cm}  ~   \\ \displaystyle 
  {{\partial \rho}\over{\partial \rho_0 }} \,&=&\, \displaystyle 
\big(   C_\theta  \,  \widetilde C_\theta  \, 
- \, \varepsilon \,   S_\theta  \,  \widetilde S_\theta  \big) \,+\, 
 \varepsilon  \,   S_\theta  \,  \widetilde S_\theta  \,\, \zeta_0    
  \\ \displaystyle   \vspace{-.3cm}  ~   \\ \displaystyle 
  {{1}\over{a}}   \, {{\partial J}\over{\partial \theta}} \,&=&\, \displaystyle 
( 1 \,-\,  \varepsilon  \,  \widetilde  \varepsilon ) \, 
  C_\theta  \,  \widetilde C_\theta  \, \,  \rho_0 
\,+\,   \big(  \varepsilon \,  \widetilde  \varepsilon \,   C_\theta  \,  \widetilde C_\theta  \, 
+ \,  \varepsilon \,   S_\theta  \,  \widetilde S_\theta  \big) \, \,   {{\sigma}\over{\sigma'}}
  \\ \displaystyle   \vspace{-.3cm}  ~   \\ \displaystyle 
  {{1}\over{a}}   \, {{\partial J}\over{\partial \rho_0}} \,&=&\, \displaystyle 
\big(   C_\theta  \,  \widetilde S_\theta  \, 
- \, \varepsilon \,   \widetilde  \varepsilon \,   S_\theta  \,  \widetilde C_\theta  \big) \,+\, 
  \varepsilon \,  \widetilde  \varepsilon \,   S_\theta  \,  \widetilde C_\theta  \, \, \zeta_0  \, .
\end{array} \right. \monendstar  
Then the  determinant $ \, \Delta \, \equiv \, {{\partial \, ( \rho , \, J)}
\over {\partial \,  (\theta ,\, \rho_0 ) }} \, \,$  can be evaluated as follows~: 

\smallskip  \noindent  $ \displaystyle
\Delta \, =    {{\partial \rho}\over{\partial \theta}} \, \, 
 {{\partial J}\over{\partial \rho_0}} \, \,-\, \, 
  {{\partial \rho}\over{\partial \rho_0 }} \, \,  
 {{\partial J}\over{\partial \theta}} \, $

\smallskip  \noindent  $ \displaystyle = \, 
 a  \,  \, \Big[  \, \Big( 
(   \widetilde  \varepsilon \, - \,   \varepsilon ) \, 
  C_\theta  \,  \widetilde S_\theta  \, \,  \rho_0 
\,+\,  \varepsilon \big(   C_\theta  \,  \widetilde S_\theta  \, 
+ \,  S_\theta  \,  \widetilde C_\theta  \big) \, \,   {{\sigma}\over{\sigma'}} \Big)
\,  \Big( \big(   C_\theta  \,  \widetilde S_\theta  \, 
- \, \varepsilon \,   \widetilde  \varepsilon \,   S_\theta  \,  \widetilde C_\theta  \big) \,+\, 
  \varepsilon \,  \widetilde  \varepsilon \,   S_\theta  \,  \widetilde C_\theta  \, \, \zeta_0  
\Big)  \, \Big] \, $ 

\smallskip  \noindent  $ \displaystyle \quad  -  \, 
 a  \,  \, \Big[  \, \Big( 
\big(   C_\theta  \,  \widetilde C_\theta  \, 
- \, \varepsilon \,   S_\theta  \,  \widetilde S_\theta  \big) \,+\, 
 \varepsilon  \,   S_\theta  \,  \widetilde S_\theta  \,\, \zeta_0 \Big) \,  
 \Big( ( 1 \,-\,  \varepsilon  \,  \widetilde  \varepsilon ) \, 
  C_\theta  \,  \widetilde C_\theta  \, \,  \rho_0 
\,+\,   \big(  \varepsilon \,  \widetilde  \varepsilon \,   C_\theta  \,  \widetilde C_\theta  \, 
+ \,  \varepsilon \,   S_\theta  \,  \widetilde S_\theta  \big) \, \,   {{\sigma}\over{\sigma'}} 
\Big)  \, \Big]  $ 

\newpage 
\smallskip  \noindent  \qquad $ \displaystyle = \, 
 a  \,  \, \Big[  \, (  \varepsilon \,  \widetilde  \varepsilon - 1) \,  C_\theta^2 \, 
\big( \widetilde C_\theta^2 \,-\,  \widetilde  \varepsilon \,  \widetilde S_\theta^2 \big) \, \rho_0 
\,+\,  \big( S_\theta^2 \,+\,   \varepsilon \, C_\theta^2 \big)
\, \big( \widetilde S_\theta^2 \,-\,  \widetilde  \varepsilon \,  \widetilde C_\theta^2 \big)
\,  {{\sigma}\over{\sigma'}} $ 

\smallskip  \noindent  \qquad $ \displaystyle \qquad  \qquad 
+  \,    \widetilde  \varepsilon \,  S_\theta^2 \, 
 \big( \widetilde C_\theta^2 \,-\,  \widetilde  \varepsilon \,  \widetilde S_\theta^2 \big)
\,  {{\sigma}\over{\sigma'}} \, \zeta_0 \, \Big] $

\smallskip  \noindent 
and due to the identity (\ref{trigo-tilde}) and the hypothesis (\ref{epsilon-tilde-carre}), 
we have the intermediate result~: 
\moneq   \label{det-jaco-2}
\Delta  \, \,=\, 
 a  \,  \, \Big[   (  \varepsilon  \,  \widetilde  \varepsilon -1 ) \, 
 C_\theta ^2 \, \, \rho_0  \,\, - \,\,  \widetilde  \varepsilon \,  
 \big(  S_\theta^2 \,+\,  \varepsilon \, C_\theta^2   \big) \,  {{\sigma}\over{\sigma'}}
  \,\, + \,\,   \widetilde  \varepsilon \,  
S_\theta ^2 \,  \,  {{\sigma}\over{\sigma'}} \, \zeta_0 \,  \, \Big] \, . 
 \monend    
Then 

\smallskip  \noindent  $ \displaystyle 
\Delta \, \,=\, \, a \,  {{\sigma}\over{\sigma'}} \,  \Big[ \,   (  \varepsilon  
\,  \widetilde  \varepsilon -1 ) \, 
 C_\theta ^2 \,   \Big(  {{\sigma^*}\over{\sigma}} \,+ \, 1 \Big) 
 \,\, - \,\,  \widetilde  \varepsilon \,  
 \big(  S_\theta^2 \,+\,  \varepsilon \, C_\theta^2   \big)
  \,\, + \,\,   \widetilde  \varepsilon \,  
S_\theta ^2 \,  \, \zeta_0 \,   \Big] $

\smallskip  \noindent  \quad $ \displaystyle \,\,
=\, a \,  {{\sigma}\over{\sigma'}} \, \Big[ \,  (\varepsilon \,  \widetilde  \varepsilon 
-1) \,\, C_\theta^2 \,\,  {{\sigma^*}\over{\sigma}} 
\,\, - \, \,(  C_\theta^2 \,+ \, \widetilde  \varepsilon \,  S_\theta^2 )
\, \, +\,\,  \widetilde  \varepsilon \, S_\theta^2 \, \zeta_0  \, \Big] $

\smallskip  \noindent
and the expression (\ref{det-jaco})  is established.  
$ \hfill \square $  
 
\smallskip  \monitem 
After determining the parameters $ \, \theta \,$ and $ \, \rho_0 \,$ 
it is possible to evaluate  the mathematical entropy $ \,\,  \eta \,\,$ thanks to the relation 
(\ref{entropie}). Then we can differentiate this mathematical entropy relatively 
the density  $ \, \rho \,$ and the momentum $ \, J . \,$ We have the following remarquable 
 coherence property.

\smallskip   \monitem   { \bf Proposition 14. Entropy variables } 

\noindent 
If the jacobian determinant  (\ref{det-jaco}) is not null, then the entropy variables 
$ \, \varphi \, \equiv \, ( \alpha , \, \beta ) \,$ 
can be determined according to  the relations (\ref{var-entrop-simple}). 
In other terms, 
\moneq   \label{var-entrop-simple-2}
\alpha \,\equiv \, {{\partial \eta}\over{\partial \rho}} \,\,=\,\, 
 \widetilde  C_\theta \,\,   \sigma'(\rho_0) \,, \qquad 
\beta  \,\equiv \, {{\partial \eta}\over{\partial J}} \,\,=\,\, 
-  \widetilde  {{\varepsilon}\over{a}}  \,\,
 \widetilde  S_\theta \,\,   \sigma'(\rho_0) \,.  
\monend

\smallskip  \noindent      {\bf Proof of Proposition 14. } 

\noindent 
We first invert the jacobian matrix~:

\smallskip  \noindent  \qquad  \qquad  \qquad  $ \displaystyle
 \begin{pmatrix}  \displaystyle  {{\partial \rho}\over{\partial \theta}} & 
 \displaystyle   {{\partial \rho}\over{\partial \rho_0}} \cr   \vspace{-.4cm}  ~ & ~  \cr  
\displaystyle  {{\partial J}\over{\partial \theta}} & 
 \displaystyle   {{\partial J }\over{\partial \rho_0}}  \end{pmatrix}^{-1}   
 \,\, = \,\, 
 \begin{pmatrix}  \displaystyle  {{\partial \theta}\over{\partial \rho}} & 
 \displaystyle   {{\partial \theta}\over{\partial J}} \cr   \vspace{-.4cm}  ~ & ~  \cr  
\displaystyle  {{\partial \rho_0}\over{\partial \rho}} & 
 \displaystyle   {{\partial \rho_0}\over{\partial J}}  \end{pmatrix} 
 \,\, = \,\, {{1}\over{\Delta}} \,\, 
 \begin{pmatrix}  \displaystyle \,\,  {{\partial J }\over{\partial \rho_0}} & 
 \displaystyle  - {{\partial \rho}\over{\partial \rho_0}} \cr   \vspace{-.4cm}  ~ & ~  \cr  
\displaystyle - {{\partial J}\over{\partial \theta}} & 
 \displaystyle \,\,  {{\partial \rho }\over{\partial \theta}} \end{pmatrix} \, .   $ 

\smallskip  \noindent 
Then 
\moneqstar    \left\{ \begin{array} [c]{rclcl} \displaystyle 
  {{\partial  \theta}\over{\partial \rho}} \,&=&\, 
\displaystyle  {{1}\over{\Delta}} \,   {{\partial  J}\over{\partial \rho_0 }}
  \,&=&\, \displaystyle   {{a}\over{\Delta}} \, \Big[ (   C_\theta  \,  \widetilde S_\theta  \,
-  \,  \varepsilon \,  \widetilde \varepsilon \,   S_\theta  \,  \widetilde C_\theta )
\,+\,  \varepsilon \,  \widetilde \varepsilon \,    S_\theta  \,  \widetilde C_\theta
\, \zeta_0 \Big] 
  \\ \displaystyle   \vspace{-.4cm}  ~   \\ \displaystyle 
  {{\partial \rho_0}\over{\partial \rho }} \,&=&\, \displaystyle 
\displaystyle  - {{1}\over{\Delta}} \,   {{\partial  J}\over{\partial \theta }}
  \,&=&\, \displaystyle  - {{a}\over{\Delta}} \,   \Big[ ( 1 -   \varepsilon \,  \widetilde \varepsilon )
\,   C_\theta  \,  \widetilde C_\theta  \, \rho_0 \,+\, 
\big(   \varepsilon \,  \widetilde \varepsilon \,  C_\theta  \,  \widetilde C_\theta  \,
+ \,   \varepsilon \,   S_\theta  \,  \widetilde S_\theta \big) \,  {{\sigma}\over{\sigma'}} 
  \Big] 
  \\ \displaystyle   \vspace{-.4cm}  ~   \\ \displaystyle 
  {{\partial \theta}\over{\partial J }} \,&=&\, \displaystyle 
\displaystyle  - {{1}\over{\Delta}} \,   {{\partial  \rho}\over{\partial \rho_0  }}
  \,&=&\, \displaystyle  - {{1}\over{\Delta}} \,    
\Big[ (  C_\theta  \,  \widetilde C_\theta  \, - \,  \varepsilon \, 
 S_\theta  \,  \widetilde S_\theta ) \, + \,  \varepsilon \, 
  S_\theta  \,  \widetilde S_\theta \, \zeta_0  \Big] 
  \\ \displaystyle   \vspace{-.4cm}  ~   \\ \displaystyle 
{{\partial \rho_0 }\over{\partial J }} \,&=&\, \displaystyle 
\displaystyle   {{1}\over{\Delta}} \,   {{\partial  \rho}\over{\partial \theta  }}
  \,&=&\, \displaystyle   {{1}\over{\Delta}} \, \Big[ (  \widetilde \varepsilon
\, - \,  \varepsilon ) \,  C_\theta  \,  \widetilde S_\theta \, \rho_0 \, + \, 
 \varepsilon  \, (   C_\theta  \,  \widetilde S_\theta \, + \, 
  S_\theta  \,  \widetilde C_\theta ) \,   {{\sigma}\over{\sigma'}} \Big] \, . 
\end{array} \right. \monendstar  
We first observe that 

\smallskip  \noindent  \qquad  \qquad    $ \displaystyle
 {{\partial \eta}\over{\partial \theta}} \,=\, 
 {{\partial}\over{\partial \theta}} \big( C_\theta \, \sigma(\rho_0) \big) \,=\,  \varepsilon 
\,   S_\theta \, \sigma \,, \qquad  
 {{\partial \eta}\over{\partial \rho_0}} \,=\, 
 {{\partial}\over{\partial \rho_0}} \big( C_\theta \, \sigma(\rho_0) \big) \,=\, 
 C_\theta \, \sigma' \, .  $

\smallskip  \noindent 
Then with  the  chain rule, we have 

\smallskip  \noindent  \qquad   $ \displaystyle
\alpha \, = \, {{\partial \eta}\over{\partial \rho}} \,=\, 
 {{\partial \eta}\over{\partial \theta}} \,  {{\partial \theta}\over{\partial \rho}} \, + \, 
 {{\partial \eta}\over{\partial \rho_0}} \,  {{\partial \rho_0}\over{\partial \rho}} $ 

\smallskip  \noindent   \qquad  $ \displaystyle \quad = \, 
  {{a}\over{\Delta}} \, \Big[ 
 \varepsilon \,   S_\theta \, \sigma \, \big(  (   C_\theta  \,  \widetilde S_\theta  \,
-  \,  \varepsilon \,  \widetilde \varepsilon \,   S_\theta  \,  \widetilde C_\theta )
\,+\,  \varepsilon \,  \widetilde \varepsilon \,    S_\theta  \,  \widetilde C_\theta
\, \zeta_0 \big) $  

 \noindent   \qquad  $ \displaystyle \qquad \qquad  
\,-\,  C_\theta \, \sigma' \, 
\big(  ( 1 -   \varepsilon \,  \widetilde \varepsilon )
\,   C_\theta  \,  \widetilde C_\theta  \, \rho_0 \,+\, 
(   \varepsilon \,  \widetilde \varepsilon \,  C_\theta  \,  \widetilde C_\theta  \,
+ \,   \varepsilon \,   S_\theta  \,  \widetilde S_\theta \,  {{\sigma}\over{\sigma'}} \big) 
\Big] $

\smallskip  \noindent  \qquad   $ \displaystyle \quad = \, 
  {{a}\over{\Delta}} \, \Big[    ( \varepsilon \,  \widetilde \varepsilon - 1 ) \, 
  C_\theta^2   \,  \widetilde C_\theta  \, \rho_0 \, \sigma' \, + \, 
   \widetilde \varepsilon \,   S_\theta^2   \,  \widetilde C_\theta  \, \zeta_0 \, \sigma 
\, - \,  \varepsilon \,  \widetilde C_\theta \, (   S_\theta^2   \, + 
 \varepsilon \,   C_\theta^2  ) \, \sigma \Big] $

\smallskip  \noindent  \qquad   $ \displaystyle \quad = \, 
  {{a}\over{\Delta}} \,  \widetilde C_\theta  \, \sigma' \, 
 \Big[    ( \varepsilon \,  \widetilde \varepsilon - 1 ) \, 
  C_\theta^2   \,  \rho_0 \, + \, 
   \widetilde \varepsilon \,   S_\theta^2   \,  \zeta_0 \,  {{\sigma}\over{\sigma'}}
\, - \,  \varepsilon \, (   S_\theta^2   \, +   \varepsilon \,   C_\theta^2  ) \, 
 {{\sigma}\over{\sigma'}}  \Big] \,\, = \,\, 
 \widetilde C_\theta  \, \sigma' \, $

\smallskip  \noindent  
 due to (\ref{det-jaco-2}). In an analogous way, we have

\smallskip  \noindent   \qquad   $ \displaystyle
\beta \, = \, {{\partial \eta}\over{\partial J }} \,=\, 
 {{\partial \eta}\over{\partial \theta}} \,  {{\partial \theta}\over{\partial J }} \, + \, 
 {{\partial \eta}\over{\partial \rho_0}} \,  {{\partial \rho_0}\over{\partial J }} $

\smallskip  \noindent   \qquad   $ \displaystyle \quad = \, 
  {{1}\over{\Delta}} \, \Big[ 
 - \varepsilon \,   S_\theta \, \sigma \, \Big( 
 (  C_\theta  \,  \widetilde C_\theta  \, - \,  \varepsilon \, 
 S_\theta  \,  \widetilde S_\theta ) \, + \,  \varepsilon \, 
  S_\theta  \,  \widetilde S_\theta \, \zeta_0  \Big)   $

 \noindent   \qquad   $ \displaystyle \qquad \qquad  
\,+\,  C_\theta \, \sigma' \, 
\Big(   (  \widetilde \varepsilon
\, - \,  \varepsilon ) \,  C_\theta  \,  \widetilde S_\theta \, \rho_0 \, + \, 
 \varepsilon  \, (   C_\theta  \,  \widetilde S_\theta \, + \, 
  S_\theta  \,  \widetilde C_\theta ) \,   {{\sigma}\over{\sigma'}} \Big) \, \Big]   $

\smallskip  \noindent  \qquad    $ \displaystyle \quad = \, 
  {{\widetilde  \varepsilon}\over{\Delta}} \, \Big[ 
( 1 -  \varepsilon \,   \widetilde \varepsilon ) \,   C_\theta^2 \,  \widetilde S_\theta \, 
\rho_0 \,  \sigma' \, \,+ \,\,   \widetilde \varepsilon
\big(  S_\theta^2 \,+\  \varepsilon \,  C_\theta^2  \big) \,  \widetilde S_\theta \, \sigma 
\,\, - \,   \widetilde \varepsilon \,  S_\theta^2  \,  \widetilde S_\theta \,
 \zeta_0  \, \sigma \Big] $ 

\smallskip  \noindent  \qquad    $ \displaystyle \quad = \, 
  {{\widetilde  \varepsilon}\over{\Delta}} \,  \widetilde S_\theta \,  \sigma' \,    \Big[ 
( 1 -  \varepsilon \,   \widetilde \varepsilon ) \,   C_\theta^2 \, 
\rho_0  \, \,+ \,\,   \widetilde \varepsilon
\big(  S_\theta^2 \,+\, \varepsilon \,  C_\theta^2  \big) \,  {{\sigma}\over{\sigma'}}
\,\, - \,   \widetilde \varepsilon \,  S_\theta^2  \,  \zeta_0  \, {{\sigma}\over{\sigma'}}  \Big]  
\,\,= \,\,  
-  {{\widetilde  \varepsilon}\over{a}} \,\,  \widetilde S_\theta \,\,  \sigma' \,    $

\smallskip  \noindent 
 due to (\ref{det-jaco-2}). The proof is completed. 
$ \hfill \square $  

\smallskip  \monitem 
We particularize our study and wish to construct  nontrivial 
hyperbolic systems  with entropy  velocity  covariant under 
the Lorentz group and - or the circular group.
It seems essential to maintain a constant sign for the determinant 
$\, \Delta \,$  presented at the relation  (\ref{det-jaco}). 
Moreover, we have to enforce the strict convexity of the
mathematical entropy $ \, ( \rho , \, J ) \longmapsto 
\eta( \rho , \, J ) .\,$ 
In this contribution, we only consider the particular case 
\moneq   \label{epsilon-epsilon-tilde}
\varepsilon \, = \,  \widetilde \varepsilon  \, .  
\monend    
We focus on ``circular elliptic'' and ``Lorentz hyperbolic''
two by two systems of conservation laws. 
Observe first that due to the condition (\ref{epsilon-epsilon-tilde}), 
the expression (\ref{variables-conservatives-3}) 
of the conserved variables simplify into 
\moneq   \label{variables-conservatives-4}
\rho \,\,=\,\, \rho_0 \,+\, \varepsilon \, S_\theta^2 \,  {{\sigma}\over{\sigma'}} 
\,, \quad 
J \,\,=\,\, a \, S_\theta \, C_\theta \,  {{\sigma}\over{\sigma'}}   \, . 
\monend    
Then the mathematical entropy $ \, \eta \,$ constrained to the condition
 (\ref{entropie}) is strictly convex if and only if the determinant 
of its hessian 
$\, {\rm det} \, (  {\rm d}^2 \eta ) \, $ 
is stricly positive and if the second derivative 
$ \, {{\partial^2 \eta}\over{\partial \rho^2}} \,$ 
is also strictly positive.

\smallskip   \monitem   { \bf Proposition 15. Technical lemma for convexity conditions  } 

\noindent  
Under the hypotheses (\ref{epsilon-tilde-carre}) and   (\ref{epsilon-epsilon-tilde}), 
the expression of the jacobian determinant $ \, \Delta \,$ is  given according to 
\moneq   \label{det-jaco-3}
\Delta  \, = \, a \, \delta \,  {{\sigma}\over{\sigma'}} \, 
\monend    
with 
\moneq   \label{petit-delta}
\delta \,=\,  \varepsilon \, S_\theta^2 \, \zeta_0 
\,\, - \, \,(  C_\theta^2 \,+ \,   \varepsilon \,  S_\theta^2 )  \, . 
\monend    
The two signed expressions are given according to  
\moneq   \label{det-hessian}
{\rm det} \, \big(  {\rm d}^2 \eta \big) \, = \, {{\varepsilon \, 
\sigma' \, \sigma''}\over{a \, \Delta}} \, \, = \, {{\varepsilon}\over{a^2}} \, 
\big( \sigma' \big)^2  \,  {{\sigma''}\over{\sigma}} \, {{1}\over{\delta}} 
\monend    
and
\moneq   \label{d2eta-sur-d2rho}
 {{\partial^2 \eta}\over{\partial \rho^2}} \,=\,
 {{\partial \alpha}\over{\partial \rho}} \,= \, {{C_\theta}\over{\delta}}
\, \Big( \varepsilon \, S_\theta^2 \, \zeta_0 \, {{\sigma'^2}\over{\sigma}}
\,- \, \big( C_\theta^2 + \varepsilon \, S_\theta^2 \big) \, \sigma'' \Big) \, . 
\monend

\smallskip  \noindent      {\bf Proof of Proposition 15. } 

\noindent 
The expression $\, \delta \,$ in (\ref{petit-delta}) is just 
a specific name for the big braket  in right hand side of 
the expression (\ref{det-jaco}) and the relation  (\ref{det-jaco-3})  
is satisfied. Under  the condition 
(\ref{epsilon-epsilon-tilde}), the expression   $\, \delta \,$ 
is clearly  reduced to the expression  (\ref{petit-delta}). 
%
%
%
%
We establish now the expression (\ref{det-hessian}) of the 
determinant of the hessian matrix. 
From the chain rule, we have     

\smallskip  \noindent  \qquad   $ \displaystyle 
{\rm det} \, \big(  {\rm d}^2 \eta \big) \, = \, 
  {{\partial \, ( \alpha , \, \beta)} \over {\partial  \, ( \rho  ,\, J )}} 
\,\,=\,\,
  {{\partial ( \alpha , \, \beta)} \over {\partial( \theta  ,\, \rho_0  )}} 
\,\, 
  {{\partial (  \theta  ,\, \rho_0 )} \over {\partial( \rho  ,\, J  )}} 
\,\, = \,\, {{1}\over{\Delta}} \, \, {\rm det} 
 \begin{pmatrix}   \displaystyle  {{\partial \alpha}\over{\partial \theta}} &
  \displaystyle  {{\partial \alpha}\over{\partial \rho_0}} \cr  \vspace{-.4cm} \cr 
  \displaystyle  {{\partial \beta}\over{\partial \theta}} & 
 \displaystyle  {{\partial \beta}\over{\partial \rho_0}}  \end{pmatrix} $

\smallskip  \noindent  \qquad    $ \displaystyle \qquad \qquad  \,\, = \,\,
  {{1}\over{\Delta}} \, \, {\rm det} 
\begin{pmatrix}   \displaystyle  \varepsilon \, S_\theta \, \sigma'  &
  \displaystyle C_\theta \, \sigma''  \cr 
  \displaystyle -  {{\varepsilon}\over{a}} \, C_\theta \, \sigma' & 
 \displaystyle  -  {{\varepsilon}\over{a}} \, S_\theta \, \sigma''   \end{pmatrix} 
 \,\, = \,\,  - {{1}\over{\Delta}} \,  {{\varepsilon}\over{a}} \,  \sigma'
\, \sigma'' \, \,  {\rm det} 
\begin{pmatrix}   \displaystyle  \varepsilon \, S_\theta    &
  \displaystyle C_\theta  \cr 
  \displaystyle \, C_\theta   & 
 \displaystyle   \, S_\theta     \end{pmatrix}  $

\smallskip  \noindent  
 due to (\ref{var-entrop-simple-2}) and  (\ref{epsilon-epsilon-tilde}). Then 
   $ \displaystyle \,\, 
{\rm det} \, \big(  {\rm d}^2 \eta \big) \, = \, {{1}\over{a \, \delta}}
\,  {{\sigma'}\over{\sigma}} \,  {{\varepsilon}\over{a}} \,  \sigma'
\, \sigma'' \,    $ 
is a consequence of the  trigonometrical property (\ref{trigo})
and  the expression (\ref{det-hessian}) follows.
The evaluation of 
(\ref{d2eta-sur-d2rho}) is also elementary.
From the chain rule, we have

\smallskip  \noindent \qquad   $ \displaystyle 
 {{\partial \alpha}\over{\partial \rho}} \,\,= \,\,
 {{\partial \alpha}\over{\partial \theta}} \, 
 {{\partial \theta}\over{\partial \rho}} \,+\, 
 {{\partial \alpha}\over{\partial \rho_0}} \, 
 {{\partial \rho_0}\over{\partial \rho}} \, \,\,= \,\,
{{a}\over{\Delta}} \, \Big[   \varepsilon \, S_\theta \, \sigma' \, S_\theta
\, C_\theta \, \zeta_0 \,-\,  C_\theta \,  \sigma'' \,
\big(     C_\theta^2 \, + \,   \varepsilon \,   S_\theta^2 \, \big)  {{\sigma}\over{\sigma'}} 
  \Big]  $

\smallskip  \noindent  \qquad    $ \displaystyle \quad  \,\, \,\,  \,\,= \,\, 
{{1}\over{\delta}} \,  {{\sigma'}\over{\sigma}} \, C_\theta 
\,  \Big[   \varepsilon \, S_\theta^2 \, \sigma' \, \zeta_0 \,-\,   \sigma'' \,
\big(     C_\theta^2 \, + \,   \varepsilon \,   S_\theta^2 \big)  \,  {{\sigma}\over{\sigma'}} 
  \Big]  $ \quad

\smallskip  \noindent 
and the expression  (\ref{d2eta-sur-d2rho}) is clear. 
$ \hfill \square $  

\smallskip   \monitem   
We detail now two explicit examples of hyperbolic system of order two.
The first one is invariant under the action of the circular group
and the second under  the Lorentz group.

\smallskip   \monitem   { \bf Proposition 16. Exponential entropy for
 a ``circular-elliptic'' system } 

\noindent  
We introduce two  strictly positive constants $ \, \rho_* \,$ 
and $ \, \overline \sigma \,$ and an exponential entropy 
on the manifold $\, \Omega_0 \,$ thanks to 
\moneq   \label{entropie-exponentielle}
\sigma (\rho_0) \,\, = \, \,    \overline \sigma \, \exp \Big( 
{{\rho_0}\over{\rho_*}} \Big) \, . 
\monend    
We set  
\moneq   \label{grand-psi}
 \left\{ \begin{array} [c]{l} \displaystyle  \displaystyle 
\Psi \,=\,   {{2\, J}\over{a \,\rho_*}}    \,,\qquad 
\rho_0  \,=\, 
 \rho + {{1}\over{2}}   \big( 1 - \sqrt{1 - \Psi^2}  \big)\, \rho_*      \,, \qquad 
\theta  \,=\,  {{1}\over{2}}     \, {\rm Arcsin} \,   \Psi   \,,
 \\  \vspace{-.4cm}    \\ \displaystyle 
u  \, = \, c \, {{J}\over{\mid \! J \! \mid}} \,\, 
\sqrt{ {{1 - \sqrt{1 - \Psi^2}}\over{1+\sqrt{1 - \Psi^2}}}}  
 \, = \, c \, {\rm sgn}(J) \, \tan \theta \,,\qquad  
 p_0  \,=\,  a \, c \, \, \big( \rho_0 \,-\, \rho_* \big)  \, .
\end{array} \right. \monend   
%
%
The system of conservation laws 
\moneq   \label{circulaire-elliptique}
 \left\{ \begin{array} [c]{rcl} \displaystyle 
{{\partial \rho}\over{\partial t}} \,+\, 
{{\partial}\over{\partial x}}  \Big(  
 {{c}\over{a}} \,  J 
\Big)   \,\, &=& \,\, 0   
 \\  \vspace{-.4cm}    \\ \displaystyle  
{{\partial J}\over{\partial t}} \,+\, 
{{\partial}\over{\partial x}} \big(  u \, J \,+ \, p_0 
\big)   \,\, &=& \,\, 0   
\end{array} \right. \monend  
is hyperbolic if we suppose 
\moneq   \label{hyp-circulaire-elliptique}
\mid \! \Psi \! \mid \,\, \leq \,\, 1 \, ,\qquad 
  \mid \! \theta \! \mid  \,\, \leq \,\, {{\pi}\over{4}} \, . 
\monend 
 It admits a mathematical entropy given by the expression 
\moneq   \label{entropie-circulaire-elliptique}
\eta ( \rho ,\, J) \,\, \equiv \, \, 
\sqrt{ {{1}\over{2}} \, \big( 1 + \sqrt{1 - \Psi^2} \big) }      \,\,\,\, 
\sigma \Big(
\rho \,+\, {{\rho_*}\over{2}} \, \big( 1 - \sqrt{1-\Psi^2}   \big) \Big) 
\monend  
and the dissipation of entropy of entropy solutions 
is given by the inequality 
\moneq   \label{ineg-entropie-circulaire-elliptique}
{{\partial \eta }\over{\partial t}} \,+\, 
{{\partial}\over{\partial x}} \big(  u \, \eta   \big)    \,\, \leq \,\, 0   \, . 
\monend

\smallskip  \noindent      {\bf Proof of Proposition 16. } 

\noindent 
With the choice (\ref{entropie-exponentielle})
of an exponential entropy, the parameter $\, \zeta_0 \,$ introduced at the relation 
(\ref{zeta_zero}) is null. Because $ \, \overline \sigma \,$
is supposed strictly positive, 
$\, \sigma \,$ is a convex function and  the condition   
$ \, {\rm det} \, \big(  {\rm d}^2 \eta \big) \,> \, 0 \,$ 
is reduced to $ \, \delta < 0 \, $ if we look to a  ``circular-elliptic'' 
hyperbolic system, {\it id est} if  we suppose $ \, \varepsilon = -1 \,$ 
(see (\ref{det-hessian})). 
Then with the notation $ \, \Psi \,$
introduced in (\ref{grand-psi}), the relations (\ref{variables-conservatives-4})
can be written 
\moneqstar   
\sin \, ( 2 \, \theta) \,\,=\,\, \Psi \, , \qquad 
\rho \,\,=\,\, \rho_0 \,-\,  \sin^2 \theta \, \rho_*  \, . 
\monendstar    
and the calculus of $ \, \theta \,$ proposed in  (\ref{grand-psi})
is natural. We have the same remark for the expression 
of  $\, \rho_0 \,$ in  (\ref{grand-psi})
as a consequence of the identity 
$\, \sin^2 \, \theta \,=\, {{1}\over{2}} \, \big( 1 - 
\sqrt{ 1 - \sin^2 (2 \, \theta)  }    \big) . \,$ 
The second condition in (\ref{hyp-circulaire-elliptique}) is a consequence
of the expression of $ \, \theta \,$ proposed in  (\ref{grand-psi}).
This last condition is also necessary when we impose
$\, \delta < 0  \, $ because the relation (\ref{petit-delta}) 
can be now written $ \,\delta = -( \cos^2  \theta - \sin^2 \theta ) .\,$  
The velocity in  (\ref{grand-psi}) 
is simply the actual formulation of the relation (\ref{vitesse-groupe})
because 
\moneqstar  
\tan \,  \theta \,=\,  
 {\rm sgn}(\theta) \, \sqrt { 
 {{1 - \sqrt{ 1 - \sin^2 (2 \, \theta)  } 
}\over { 1 + \sqrt{ 1 + \sin^2 (2 \, \theta)  } }}} . \,  
\monendstar
%
%

 \monitem
The static pressure $ \, p_0 \,$ is   related to the entropy $ \, \sigma \,$
thanks to (\ref{pression}). 
The expressions (\ref{def-flux-thermo}) and  (\ref{thermo-flux-2}) 
of the physical and thermodynamix fluxes, 
the algebraic forms (\ref{matrice-Y-theta})  of the matrix $ \, Y_\theta \,$
and the choice (\ref{g-zero-Neq2}) for $ \, g_0 \,$ show that 
the first component $ \, f_1  \, $ of $ \, f(W) \,$ is given by

\smallskip  \noindent  \qquad   $ \displaystyle 
f_1 \,\,=\,\, \rho \, u \, - {{\sin   \theta}\over{a\, \cos \, \theta }} \, p_0    
 \,\,=\,\,  \rho \, u \, -  {{1}\over{a \, c}} \, u \, a \, c \, ( \rho_0 - \rho_* ) 
 \,\,=\,\,  \rho_* \, u \, (1 - \sin^2 \theta)    $

\smallskip  \noindent   \qquad     $ \displaystyle \quad \,\,\,  = \, 
\rho_* \, c \, \sin   \theta  \, \cos  \theta
 \,\,=\,\, {{1}\over{2}} \,  \rho_* \, c \, \sin (2 \, \theta) 
 \,\,=\,\, {{1}\over{2}} \,   \rho_* \, c \, \Psi \, 
 \,\,=\,\, {{1}\over{2}} \,   \rho_* \, c \,  {{2\, J}\over{a \,\rho_*}}  
 \,\,=\,\, {{c}\over{a}} \,  J \, .  $

\smallskip  \noindent 
Then the relations (\ref{circulaire-elliptique}) are established.
The expression  (\ref{entropie-circulaire-elliptique}) of the mathematical 
entropy is then a consequence of (\ref{entropie})~:
$ \, \eta = \cos \theta \, \, \sigma(\rho_0) .\,$ 
A particular hyperbolic system   with entropy velocity 
 invariant under  the action of  the circular group is constructed. 
$ \hfill \square $  

\smallskip  \newpage  \monitem   { \bf Proposition 17. Homographic  entropy for
 a ``Lorentz-hyperbolic'' system } 

\noindent  
We consider  two  strictly positive constants $ \, \rho_* \,$ 
and $ \, \overline \sigma \,$ and an homographic  entropy $\, \sigma \,$ 
on the manifold $\, \Omega_0 \,$ defined  by the relation  
\moneq   \label{entropie-homographic}
\sigma (\rho_0) \,\, = \, \,    - \, \overline \sigma \, 
{{\rho_0}\over{\rho_0 + \rho_*}}  \,,  \qquad \rho_0 > 0 \, . 
\monend    
We set  
\moneq   \label{grand-phi-hyperbolic} 
\Phi (\rho_0) \,\, \equiv \,\,   {{2\, J}\over{a}}  \, \, 
 {{\sigma'}\over{\sigma}} 
\,\, =  \,\,  {{2\, J}\over{a}}  \,  {{\rho_*}\over{\rho_0 \, ( \rho_0 + \rho_*)}}   \,  
 \monend   
and   we suppose 
\moneq   \label{hyp-lorentz-hyperbolic}
{{\mid \! J  \! \mid}\over{a}} \,\, \leq \,\, 
\rho \,\, \leq \,\, {{1}\over{2}} \, \rho_*   \, . 
\monend 
Then the equation 
\moneq   \label{lorentz-equation}
F(\rho_0) \, \equiv \,  \rho_0 +  {{\sigma}\over{2 \, \sigma'}} 
\Big( \sqrt{1 + \Phi (\rho_0)^2} - 1 \Big) \,\,=\,\, \rho 
\monend 
has a unique solution $ \, \rho_0 \,$ such that 
\moneq   \label{condition-rho-zero}
0 \,\leq \,  \rho_0 \,\leq \,  {{1}\over{2}} \, \rho_* \, . 
\monend 
We consider the expressions 
\moneq   \label{expressions-hyperbolic} 
p_0 \,=\, a \, c \, {{\rho_0^2}\over{\rho_*}} \,, \quad 
\theta \,=\, {{1}\over{2}} \, {\rm Argsh} \,  \Phi(\rho_0)  \,, \quad 
u  \,=\, a \, \tanh \theta \,=\, c \,  {\rm sgn}(J) \, 
\sqrt{ {{ \sqrt{1 + \Phi^2} - 1 }\over{ \sqrt{1 + \Phi^2} + 1 }}  } \, . 
\monend 
Then the system of conservation laws 
\moneq   \label{lorentz-hyperbolic}
 \left\{ \begin{array} [c]{rcl} \displaystyle 
{{\partial \rho}\over{\partial t}} \,+\, 
{{\partial}\over{\partial x}}  
\Big[ \Big(   \rho + {{\rho_0}\over{\rho_*}} \Big) \, u \Big] 
  \,\, &=& \,\, 0   
 \\  \vspace{-.4cm}    \\ \displaystyle  
{{\partial J}\over{\partial t}} \,+\, 
{{\partial}\over{\partial x}} \big(  u \, J \,+ \, p_0 
\big)   \,\, &=& \,\, 0   
\end{array} \right. \monend  
is hyperbolic and   admits a mathematical entropy given by the expression 
\moneq   \label{entropie-lorentz-hyperbolic}
\eta ( \rho ,\, J) \,\, \equiv \, \, 
\sqrt{ {{1}\over{2}} \, \big( 1 + \sqrt{1 + \Phi(\rho)^2} \big) }    
  \,\,\,  \sigma ( \rho_0 )  \, . 
\monend  
The dissipation of entropy of entropy solutions 
is still given by the inequality (\ref{ineg-entropie-circulaire-elliptique}). 
%

\smallskip  \noindent      {\bf Proof of Proposition 17. } 

\noindent 
In this case of the Lorentz group, we have $\, \varepsilon = 1 , \,$
$ \, S_\theta \equiv \sinh \theta \,$ and 
$ \, C_\theta \equiv \cosh \theta . \,$ 
We have chosen with (\ref{entropie-homographic})
a negative convex function $ \, \sigma .  \,$ 
Then ({\it c.f.} (\ref{det-hessian}))  the condition 
$ \, \delta < 0 \,$  has to be satisfied to assume 
$ \, {\rm det} \, \big(  {\rm d}^2 \eta \big) \,> \, 0 .\,$ 
We have  
\moneqstar
  {{\sigma}\over{\sigma'}}  \,=\, {{\rho_0}\over{\rho_*}} \, 
( \rho_0 + \rho_* ) \,, \quad 
\zeta_0 \,\equiv\, {{\partial}\over{\partial \rho_0}}
\Big(  {{\sigma}\over{\sigma'}} \Big) 
\,=\, 1 \,+\, 2 \,  {{\rho_0}\over{\rho_*}} \, . 
\monendstar    
The relations (\ref{variables-conservatives-4}) can now be written as 
\moneqstar   
\sinh  \, ( 2 \, \theta) \,\,=\,\, \Phi (\rho_0) \, , \qquad 
\rho \,\,=\,\, \rho_0 \,+\,  \sinh^2 \theta \, \,  {{\sigma}\over{\sigma'}}  \, . 
\monendstar    
The identity $ \, \sinh^2 \theta \,\equiv \, {{1}\over{2}} \big( 
\sqrt{1 + \sinh^2 (2 \, \theta) } - 1 \big) \, $ 
shows that the density $ \, \rho_0 \,$ at null velocity
associated with an arbitrary state $ \, W \equiv (\rho , \, J)^{\rm \displaystyle t} \,$
has to solve the equation
\moneqstar 
   \rho \,-\, \rho_0  
\,=\,  {{1}\over{2}} \Big( \sqrt{1 + \Phi (\rho_0)^2} - 1 \Big)
\, \,   {{\sigma}\over{\sigma'}} 
\monendstar    
which is equivalent to  (\ref{lorentz-equation}). 
We remark that $ \, F(0) =  \, \mid \! J \! \mid \!  / a  \,$ 
and the first condition of (\ref{condition-rho-zero}) is natural. 

 \monitem  
We prove now that the equation  (\ref{lorentz-equation})
has a unique solution if the conditions (\ref{hyp-lorentz-hyperbolic}) 
are satisfied.  We have 

\smallskip  \noindent    $ \displaystyle  \qquad 
 {{\partial}\over{\partial \rho_0}} \Big[  {{\sigma}\over{2 \, \sigma'}}
\, \sqrt{ 1 + \Phi^2} \Big] \,\,= \,\,
 {{1}\over{2}}  {{\partial}\over{\partial \rho_0}} 
\sqrt { \Big(  {{\sigma}\over{\sigma'}} \Big)^2 
\,+\, {{4 \, J^2}\over{a^2}} } 
 \,\,= \,\, {{2}\over{4 \, \sqrt { \Big(  \displaystyle  {{\sigma}\over{\sigma'}} \Big)^2 
\,+\, {{4 \, J^2}\over{a^2}} \,  } }} \, \,   {{\sigma}\over{\sigma'}}  \,
\, \zeta_0  \,\, \,\,  > \,\, 0 $

\smallskip  \noindent 
if $ \, \rho_0 > 0 . \,$ We have also 
   $ \displaystyle \quad 
 {{\partial}\over{\partial \rho_0}} \Big( \rho_0 \,-\, 
 {{\sigma}\over{2 \, \sigma'}} \Big) 
\,\,= \,\,  1 - {{\zeta_0}\over{2}} 
\,\,= \,\,  {1\over2} - {{\rho_0}\over{\rho_*}}   \,\, > \,\, 0 \quad $ 
under the condition (\ref{condition-rho-zero})  for $\, \rho_0 .\, $ 
Then we have $ \, F'(\rho_0) > 0 \,$ and the equation (\ref{lorentz-equation})
has a unique solution when  (\ref{condition-rho-zero}) is satisfied. 
The maximum of the function $ \, F(\smb) \,$ on the interval $ \, [ 
0 ,\, {{1}\over{2}} \, \rho_* ] \,$ can be estimated as follows. We have first
$ \quad  {{\sigma}\over{\sigma'}} \big(  {{\rho_*}\over{2}} \big) \,=\, 
 {{1}\over{2}} \, \big(  {{1}\over{2}} + 1 \big) \, \rho_* \, = \, 
 {{3}\over{4}} \, \rho_*  . \, $  Then 

\smallskip  \noindent    $ \displaystyle  \qquad 
F \Big(  {{\rho_*}\over{2}} \Big) \,\,= \,\, 
  {{\rho_*}\over{2}} \,-\,  {{3}\over{8}} \, \rho_* 
 \,+\,   {{\rho_*}\over{2}} \,
\sqrt{ \Big(  {{3}\over{4}}   \Big)^2 \,+\, \Big(  {{2 \, J}\over{a \, \rho_*}}  \Big)^2 \, }  
\,\,\, \geq  \,  \,\,   {{\rho_*}\over{2}} \,-\,  {{3}\over{8}} \, \rho_* 
\,+\,  {{3}\over{8}} \, \rho_*    \,\,= \,\,   {{\rho_*}\over{2}}   
$

\smallskip  \noindent 
and the second condition of (\ref{hyp-lorentz-hyperbolic}) is natural. 
The pressure $ \, p_0 \,$ is given by the relation  (\ref{pression})~:

\smallskip  \noindent    $ \displaystyle   \qquad  \qquad 
p_0 \,=\, - \, a \, c \,  {{\sigma^*}\over{\sigma'}}  \,=\, 
a \, c\, \Big(  {{\sigma}\over{\sigma'}} - \rho_0 \Big)  \,=\, 
a \, c\, {{\rho_0^2}\over{\rho_*}}   $

\smallskip  \noindent  
and the first relation of (\ref{expressions-hyperbolic})   is established. 
The two other relations of  (\ref{expressions-hyperbolic})
are an elementary consequence of the previous considerations.

 \monitem 
We now establish the algebraic form  of the hyperbolic system (\ref{edp}). 
From  (\ref{def-flux-thermo}) and  (\ref{thermo-flux-2}), 
the particular expression  (\ref{matrice-Y-theta})  of the matrix $ \, Y_\theta \,$
and the choice $ \, g_0 \equiv ( 0 , \, p_0 )^{\rm \displaystyle t}  \,$
proposed in  (\ref{g-zero-Neq2}) 
show that  the first component $ \, f_1  \, $ of $ \, f(W) \,$ is given by

\smallskip  \noindent    $ \displaystyle \qquad  \qquad 
f_1 \,\,=\,\, \rho \, u \, + {{\sinh   \theta}\over{a\, \cosh \, \theta }} \, p_0    
 \,\,=\,\,  \rho \, u \, +  {{1}\over{a \, c}} \,  \,u  \,\, a \, c  \, \,  {{\rho_0^2}\over{\rho_*}}
 \,\,=\,\,  \Big(  \rho \, + \,  {{\rho_0^2}\over{\rho_*}} \Big) \, u  \, .  $
 
\smallskip  \noindent 
Then  the fist line of the relations (\ref{lorentz-hyperbolic}) is established. 
The second line is straightforword to  evidence. 
The expression of the mathematical entropy 
is a direct consequence of (\ref{entropie}) and the explitation
of $\, \cosh \theta \,$ as a function of $ \, \Phi(\rho_0) \, $
as proposed in (\ref{expressions-hyperbolic}). 
This remark  completes   the construction of an  hyperbolic system   with entropy velocity 
(\ref{lorentz-hyperbolic})  invariant under  the action of  the Lorentz  group.  
$ \hfill \square $  

\bigskip \bigskip   \noindent {\bf \large 7) \quad  Conclusion}  
 

\monitem 
In this contribution, we have extended our previous work \cite{Du01}
relative to the construction of hyperbolic systems covariant under the Galileo group
and space reflection. 
We have precised the notion of ``entropy  velocity'' as an intrinsic velocity 
associated with the conservation of entropy for regular solutions. 
In particular, the usual velocity of gas dynamics system
is an  entropy  velocity. 
An  hyperbolic system of conservation laws and the 
associated conservation of entropy are covariant under the action of
a group of space-time transformations if the algebraic form
of the equations is not changed after the action of the corresponding group. 
We have proposed sufficient conditions
in order to satisfy a  group covariance property. 
These conditions introduce naturally a representation of the group
in the space of states that can be chosen in multiple ways. 
 We have also introduced a natural manifold of ``null velocity''. 
When thermodynamic data are given on this manifold, the natural question 
is the construction of the whole  hyperbolic system from the covariance property. 
We have also  natural constraints for the entropy variables. 

\monitem 
After the study of  the particular case 
of Galileo covariance  in \cite{Du01},   
we have   particularized in this contribution  our study 
for the two-dimensional case with the Lorentz and circular groups. 
We have proposed  nontrivial hyperbolic systems  covariant under
each of these two groups.  
%
For the circular example, we find 
a condition  (\ref{hyp-circulaire-elliptique}) that shows that the associated hyperbolic 
system is not covariant under the action of the whole circular group. 
The question is set \cite{Va13} to know if there exists  a mathematical obstruction
to find hyperbolic systems covariant under the action of 
the whole circular group. 
%
A better knowledge of such  systems of conservation laws 
is a possible  next  step. Last but not least,  the link with 
relativistic gas dynamics (see {\it e.g.} Souriau  \cite{So77, So78}, 
Chiu \cite{Ch73} 
and Vall\'ee \cite{Va87}) is also a natural question. 
The extension of this work to several space dimensions 
is, in principle, straightforward.

%

\bigskip  \bigskip    \noindent {\bf \large Acknowledgments}   


 \noindent    
The author  thanks  G\'ery de Saxc\'e and Claude Vall\'ee for their kind invitation to
participate to the 56th  and 57th  
CITV in august 2012 and august 2013, 
and for stimulating discussions during these  workshops.
In november 2014, Claude Vall\'ee passed away. 
Of course, this contribution is also dedicated to him. 

\bigskip  \bigskip   \noindent {\bf \large  References } 

\medskip

\end{document}